\documentclass{amsart}
\usepackage{a4wide}

\usepackage{textcmds} 
\usepackage{amsmath, amssymb, amsfonts, amstext, amsthm, amscd, mathrsfs, mathscinet}
\usepackage{mathtools}
\usepackage{verbatim}
\usepackage{graphicx}
\usepackage{color}
\usepackage{pinlabel}
\usepackage{mathrsfs} 
\usepackage{xypic}
\usepackage{enumitem}
\usepackage{tikz}
\usetikzlibrary{matrix}

\usepackage{etoolbox}
\makeatletter
\let\ams@starttoc\@starttoc
\makeatother
\usepackage[parfill]{parskip}
\makeatletter
\let\@starttoc\ams@starttoc
\patchcmd{\@starttoc}{\makeatletter}{\makeatletter\parskip\z@}{}{}
\makeatother

\usepackage{hyperref}

\newcommand{\CC}{\mathbf{C}}
\newcommand{\C}{\mathbb{C}}
\newcommand{\Z}{\mathbb{Z}}
\newcommand{\R}{\mathbb{R}}

\newcommand{\del}{\partial}
\newcommand{\delbar}{\overline{\partial}}

\newcommand{\OP}{\operatorname}

\renewcommand{\Re}[1]{\mathfrak{Re}\,#1}
\renewcommand{\Im}[1]{\mathfrak{Im}\,#1}
\renewcommand{\Re}{\mathfrak{Re}}
\renewcommand{\Im}{\mathfrak{Im}}

\newsavebox{\textvisiblespacebox}
\savebox{\textvisiblespacebox}{\texttt{aa}}
\newcommand\vartextvisiblespace[1][\wd\textvisiblespacebox]{%
  \makebox[#1]{\kern.1em\rule{.4pt}{.3ex}%
  \hrulefill%
  \rule{.4pt}{.3ex}\kern.1em}%
}

\numberwithin{equation}{section}

\newtheorem{thm}{Theorem}[section]
\newtheorem{lma}[thm]{Lemma}
\newtheorem{prp}[thm]{Proposition}
\newtheorem{cor}[thm]{Corollary}

\newtheoremstyle{TheoremNum}
    {\topsep}{\topsep}              
    {\itshape}                      
    {}                              
    {\bfseries}                     
    {.}                             
    { }                             
    {\thmname{#1}\thmnote{ \bfseries #3}}
\theoremstyle{TheoremNum}

\theoremstyle{definition}
\newtheorem{dfn}[thm]{Definition}

\newtheorem{ex}[thm]{Example}

\theoremstyle{remark}
\newtheorem{rmk}[thm]{Remark}

\begingroup 
\makeatletter 
\@for\theoremstyle:=definition,remark,plain,TheoremNum\do{%
\expandafter\g@addto@macro\csname th@\theoremstyle\endcsname{%
\addtolength\thm@preskip\parskip 
}%
} 
\endgroup 

\title[Rationally Convex Surfaces with hyperbolic tangencies]{Rationally Convex Surfaces with hyperbolic complex tangencies}

\author{Georgios Dimitroglou Rizell}
\address{Department of Mathematics\\
Uppsala University\\
Box 480\\
SE-751 06 UPPSALA\\
SWEDEN}
\email{georgios.dimitroglou@math.uu.se}
\author{Mark G. Lawrence}
\email{mgl\_61@yahoo.com}
\thanks{The first author is supported by the Knut and Alice Wallenberg Foundation through the grants KAW 2021.0191 and KAW 2023.0294.}

\begin{document}

\begin{abstract}
We construct the first examples of rationally convex surfaces in the complex plane with hyperbolic complex tangencies. In fact, we give two very different types of rationally convex surfaces: those that admit analytic fillings by handle-bodies, and those that do not have any compact Riemann surfaces attached at all. The fillable examples all live in the round sphere and are unknotted, while the non-fillable examples can moreover be produced in several different smooth isotopy classes.
\end{abstract}

\maketitle
\setcounter{tocdepth}{1}
\tableofcontents

\section{Introduction and results}

 This work concerns the construction of rational convex surfaces in $\C^2$ of different types, by the use of techniques from symplectic topology. Recall that a compact subset $K \subset \C^n$ is \textbf{rationally convex} if there exists a proper analytic variety passing through any point in the complement. It satisfies the stronger notion of being \textbf{polynomially convex} if, for any point in the complement $\OP{pt} \in \C^n\setminus K$, there exists a polynomial $p$ whose value satisfies $|p(\OP{pt})| > \max_{\mathbf{z} \in K}|p(\mathbf{z})|$. Any compact Riemann surface with boundary contained in $K$ is thus a subset of its polynomial hull -- the smallest polynomially convex subset containing the original subset $K$. In addition, existence of certain nice families of holomorphic curves, so-called holomorphic fillings, has been shown to imply rational convexity; see work by Duval--Gayet \cite{Duval:Riemann}. For this reason, understanding the Riemann surfaces that have boundary contained in $K$ is very central for the study of both rational and polynomial convexity. 

The theory of Lagrangian tori in symplectic manifolds holds a distinguished place in symplectic topology. In particular, their study has motivated a lot of research in low dimensional topology, and has also provided connections to algebraic geometry through homological mirror symmetry. A key question is deciding whether there is a filling by holomorphic discs. From the point of view of complex analysis, these questions hold interest as part of the theory of polynomial hulls of submanifolds of $\mathbb{C}^n$. 

The most central connection between symplectic topology and analytic function theory that we are interested in here is the following key fact: a Lagrangian torus for a global K\"{a}hler form on $\C^n$ is rationally convex. Conversely, if a torus is isotropic for some global K\"{a}hler form, then it is rationally convex. Both results are contained in Duval--Sibony's seminal work \cite{DuvalSibony}. Recall that the property of being Lagrangian for some K\"{a}hler form implies being totally real, but that the converse does not necessarily hold.

Suppose now that $\Sigma^2\subseteq \mathbb{C}^2$ is an embedded surface. Consider a vector-field $V$ that is tangent to $\Sigma$ and which has an algebraic number $2-2g$ of zeroes contained in the totally-real locus of $\Sigma$. If we push off $\Sigma$ along $JV$ to $\Sigma_{JV}$, we thus get a contribution of $2g-2$ from the zeros of $V$ to the algebraic intersection number $\Sigma \bullet \Sigma_{JV}=0$. If $JT\Sigma \cap T\Sigma$ is generic, it follows that there must be additional intersections that give rise to an algebraic count $2-2g$. These intersections must all come from the presence of complex tangencies; the hyperbolic (resp. elliptic) complex tangents are those generic and isolated complex tangencies that contribute to $-1$ (resp. $+1$) to this intersection number. If there are elliptic complex tangents, then there must exist a local holomorphic filling, which precludes rational convexity; see e.g.~\cite{BedfordKlingenberg}. On the other hand, if $\Sigma^2$ has only hyperbolic complex tangents, which thus necessarily means that $2g-2 \ge 0$, then the surface could be rationally convex. In Subsection \ref{sec:model} we show that such a surface can be assumed to be Lagrangian for a K\"{a}hler form that becomes degenerate precisely at the hyperbolic complex tangencies (no such K\"{a}hler form can exist when there are elliptic complex tangencies). In other words, these closed orientable high-genus surfaces can be seen as a sort of analogue of a high-genus Lagrangian surface in $\C^2$ (which cannot exist for a symplectic form that is everywhere non-degenerate by the above computation of intersection numbers).

We give two very different constructions of rationally convex surfaces in $\C^2$ with only hyperbolic points. The first construction described in Section \ref{sec:fillable} is very rigid, e.g.~since it bounds an embedded handle-body and thus lives in a unique smooth isotopy class when considered in $\C^2$. More precisely, in Subsection \ref{proof:ratconvexS3} we prove the following.
\begin{thm}
\label{thm:ratconvexS3}
There exist rationally convex oriented surfaces of genus $g \ge 1$ contained in $S^3 \subset \C^2$, whose complex tangencies consist of a number $2(g-1)$ of hyperbolic points. Furthermore, for each $g \ge 1$, such surfaces can be realised in several smooth isotopy classes of embeddings in $S^3$. This includes the class of the genus-$g$ surface in $S^3$ that gives its standard Heegaard-splitting into a pair of handle-bodies.
\end{thm}

The examples are constructed by starting with a totally real rationally convex torus in $S^3$ that admits a filling in $D^4$ by holomorphic discs. We then perform a surgery along Legendrian arcs in $S^3$ in order to increase the genus, thus adding pairs of hyperbolic complex tangencies for each handle; see Section \ref{sec:clifford} for the case of a surgery on the Clifford torus, which produces the surface in the Heegaard splitting of $S^3$, and see Section \ref{sec:knotted} for the case of a torus that is knotted inside $S^3$. (Note that a torus inside $S^3$ bounds a handle-body and is thus never knotted when considered inside $\C^2 \supset S^3$.) Even though these are high genus surfaces, they are obtained from Lagrangian tori, and consequently inherit many of the rigidity properties of Lagrangian tori. In particular, they admit plenty of holomorphic Maslov-2 discs, and they are smoothly unknotted when considered in $\C^2$ . Rationally convex tori that are totally real are Lagrangian for a global K\"{a}hler form; Lagrangian tori are smoothly unknotted by a result of first author joint with Goodman--Ivrii \cite{Dimitroglou:Isotopy}.

\begin{ex}
The surgery construction can also be used to produce non-rationally convex surfaces in $S^3$ with only hyperbolic complex tangencies. Consider e.g.~the totally real torus of vanishing Maslov class that was exhibited in \cite{Duval:Riemann}. This torus is not rationally convex since it admits a holomorphic annulus with null-homologous boundary inside the torus. Perform a surgery as described in Section \ref{sec:fillable} in order to add handles to the torus that contain hyperbolic complex tangencies. Adding this handle in a suitable position, it will neither affect the boundary of the annulus, nor its homology class. Thus we can produce genus-$g$ surfaces in $S^3$ with $g>1$ that have only hyperbolic complex tangencies, and which are not rationally convex for the same reason as the original torus; namely, they admit holomorphic annuli whose boundaries are contained inside the surface, such that the boundaries moreover are null homologous inside the surface.
\end{ex}

The second construction of rationally convex surfaces $\Sigma \subset \C^2$ with only hyperbolic complex tangencies is given in Section \ref{sec:nonfillable} has very different properties compared to the first one. Notably, these surfaces can be realised in different smooth isotopy classes in $\C^2$, and they do not admit any analytic filling.

We say that a compact Riemann surface with boundary is \textbf{attached to $\Sigma$} if there is a continuous map from the Riemann surface that takes the boundary into $\Sigma$, such that the map is holomorphic in the interior. We show that:
\begin{thm}[Theorem \ref{thm:exactsurfaces}, Corollary \ref{cor:rationallyconvex}, and Theorem \ref{thm:knotted}]
For any $g \ge 2$ there exist rationally convex smooth embeddings $\Sigma_g \subset \C^2$ of a genus-$g$ surface with $2(g-1)$ hyperbolic complex tangencies, such that there are no non-constant compact Riemann surfaces attached to $\Sigma_g$. Furthermore, for $g \gg 0$ sufficiently large, these surfaces can be realised in arbitrarily many different smooth isotopy classes.
\end{thm}
The reason for the non-existence of non-constant Riemann surfaces attached to these rationally convex surfaces is that they are exact Lagrangians for a K\"{a}hler form that is degenerate precisely at the hyperbolic points. Gromov's result \cite{Gromov:Pseudo} implies that there are no exact Lagrangians inside $\C^2$ for a non-degenerate K\"{a}hler form. However, our rationally convex surfaces are obtained from the \emph{singular} exact Lagrangians of genus $g \ge 2$ produced by Lin in \cite{Lin}.

\section{Fillable rationally convex surfaces with only hyperbolic tangencies}

\label{sec:fillable}
Here we describe a general type of surgery that can be performed to surfaces in contact three-manifolds that add a handle containing a pair of hyperbolic complex tangencies to the surface; see Lemma \ref{lma:tangencies} below. We will focus on the case when the contact manifold is the standard $(S^3,\xi)$, but the construction works in general contact three-manifolds. The case of interest to us is when the contact distribution is a field of complex hyperplanes and the original surface is totally real, i.e.~has no tangencies to the contact distribution. In this case, our surgery adds handles containing pairs of hyperbolic complex tangencies. 

Assume that we are given a smoothly embedded surface $\Sigma^2 \subset (S^3,\xi)$. For any pair of disjoint points $\{p_{ -1},p_{ 1}\} \subset \Sigma$ on the surface that satisfy the property that $T_{p_i}\Sigma \neq \xi$ for $i=\pm1$ (i.e.~$\Sigma$ is totally real near $p_i$) one can find smooth Legendrian arcs $\gamma \subset (S^3,\xi)$ with
\begin{itemize}
\item boundary at the two points $\partial \gamma=\{p_{-1},p_1\}$, where the arc moreover is transverse to $\Sigma$; and
\item interior disjoint from $\Sigma$.
\end{itemize}
Recall that there exist plenty of Legendrians inside any contact manifold; see e.g.~\cite{Etnyre:Legendrian} where it is shown that one can $C^0$-approximate any smooth knot by a Legendrian knot in the same smooth isotopy class.

Given a choice of Legendrian arc $(\gamma,\partial \gamma) \subset (S^3,\Lambda)$ as above connecting the two points $p_{-1},p_1 \in \Sigma$, we construct a standard neighborhood of the arc $\gamma$ in the following manner. First note that, since $\Sigma$ is totally real near $p_i$, its characteristic distribution $\xi \cap T\Sigma$ locally defines a foliation near $p_{\pm1}$. Recall Giroux's Reconstruction Result \cite[Proposition II.1.2(b)]{Giroux1991}, by which the characteristic distribution on a surface determines its contactomorphism class uniquely. In other words, we can identify neighbourhoods of $p_{\pm1} \in (S^3,\xi)$ with neighborhoods of the points $(\pm 1, 0,0) \in \left(\R^3_{xyz},\ker (dz-ydx)\right)$ so that the surface moreover is identified with open subsets of the planes 
$$ P_{\pm1} \coloneqq \{x=\pm 1\} \subset \left(\R^3_{xyz},\ker (dz-ydx)\right).$$
If we moreover deform the arc $\gamma$, so that it coincides with the Legendrian arc
$$[-1,1]_x \times \{y=0\} \times \{z=0\} \subset \left(\R^3_{xyz},\ker (dz-ydx)\right),$$
near $P_{\pm1}$, then one can readily use the standard neighbourhood theorem for Legendrians \cite{Geiges:Intro} to create a contactomorphic identification of a neighbourhood of (the deformed) arc $\gamma$ with a neighbourhood $U$ of
$$[-1,1]_x \times \{y=0\} \times \{z=0\} \subset \left(\R^3_{xyz},\ker (dz-ydx)\right),$$
so that $\gamma$ becomes identified with $[-1,1] \times \{(0,0)\}$, $p_{\pm1}$ becomes identified with $(\pm 1,0,0),$ and neighbourhoods of $p_{\pm1} \in \Sigma$ become identified with neighborhoods $\Sigma_{\pm 1} \cap U$ of $(\pm 1,0,0) \subset P_{\pm1}=\{x=\pm 1\}$.

Denote by $\Sigma^U \coloneqq U \cap \left( \Sigma_{-1} \cup \Sigma_{+1} \right)$ and note that these planes are pre-Lagrangian, and that their characteristic distribution is given by $\R\partial_y$. Recall that a submanifold is \textbf{pre-Lagrangian} if it admits a Lagrangian section inside the symplectisation. For a surface inside a contact three-manifold, this is equivalent to the characteristic distribution being given as the kernel of a globally defined closed one form. 

We then construct an embedded cylinder of the form
$$ C_\epsilon \coloneqq \left\{ \frac{1}{\epsilon}(y^2+z^2)=\rho(x); \: x\in[-1,1] \right\} \subset \left(\R^3,dz-ydx\right)$$
where $\rho \colon [-1,1] \to (0,1]$ is smooth inside $(-1,1)$ and satisfies
\begin{itemize}
\item $\rho(\pm x)=\rho(x)$,
\item $\rho(\pm1)=1$, and
\item $\rho''(x) > 0$ for all $x \in (-1,1)$.
\end{itemize}
In addition, we assume that $\rho(x)$ has been chosen to satisfy the additional property that
$$\Sigma^U_\epsilon \coloneqq C_\epsilon \cup \left(\Sigma^U \setminus \left\{x \in \{\pm 1\},\:\frac{1}{\epsilon}(y^2+z^2) \le \rho(x)\right\}\right)$$
is a smoothly embedded cylinder for any $\epsilon>0$ sufficiently small. (Roughly speaking, this means that $\rho'(x)$ has to blow up sufficiently fast as $x \to \pm 1$.) In particular, this cylinder coincides with $\Sigma^U$ near the boundary of the neighbourhood $U$. It will be crucial to understand the characteristic distribution on $\Sigma^U_\epsilon$ in a neighborhood of the cylinder; it is depicted in Figure \ref{fig:sigma1}.

 Given any surface $\Sigma$ as above, replacing the neighborhood $\Sigma^U$ of the two points $p_{\pm 1}$ above with the cylinder $\Sigma^U_\epsilon$ constructed using the Legendrian arc $\gamma$ yields a new smoothly embedded surface $\Sigma_\epsilon$ in the contact manifold. This surface is obtained by ``ambient 0-surgery" (or self-connected sum) performed on $\Sigma$ inside the contact manifold on the embedded 0-sphere $\{p_1,p_{-1}\} \subset \Sigma$. We call this construction an \textbf{ambient 0-surgery along the Legendrian arc $\gamma$ with parameter $\epsilon>0$}. Next, we proceed to investigate the characteristic distribution on the deformed part of the surface, i.e.~the cylinder $\Sigma^U_\epsilon$, that the construction gives rise to.

 One can readily verify that the only complex tangencies in $\Sigma^U_\epsilon$ all are contained inside the cylinder $C_\epsilon \subset \Sigma^U_\epsilon$. Moreover, the complex tangencies consist of the two hyperbolic points
$$ h_\pm \coloneqq (x,y,z)=(0,0,\pm \sqrt{\epsilon \rho(0)}) \in C_\epsilon,$$
relative the above coordinates. More details will be given in Lemma \ref{lma:tangencies} below.

\subsection{Orientation conventions}

 Here we give a detailed description how the characteristic distribution is oriented, depending on a choice of orientation of the surface and a choice of co-orientation of the contact distribution. 

Recall that $\xi=\ker\alpha$ has an orientation induced by the exterior derivative $d\alpha$ for a \emph{choice} of contact form $\alpha$ , and that we endow $Y=S^3$ with the standard orientation for which $\alpha \wedge d\alpha$ is a positive volume form. In other words, the contact manifold is oriented by $TY=\langle \OP{Reeb} \rangle \oplus \xi$, where $\OP{Reeb}$ is the non-zero Reeb vector field which is in the kernel of $d\alpha$ and normalised by $\alpha(\OP{Reeb})=1$, and $\xi$ is endowed with the orientation induced by $d\alpha$. Since we are in dimension three, the orientation of $S^3$ is independent of the choice of co-orientation of $\xi,$ i.e.~sign of $\alpha$. However, the orientation of $\xi$ depends on this choice. Here we will always consider the co-orientation induced by standard round contact form $\alpha$ on $S^3$, which induces the complex orientation on the contact planes $\xi$, and the standard orientation of $S^3$, coinciding with the boundary orientation of the complex ball. 

 For convenience, when considering the contact 3-manifold $(Y,\xi)$, we will fix a choice of contact form $\alpha$ as well. Given a choice of orientation of a surface $ S \subset (Y,\xi)$ we get an induced orientation of the \textbf{characteristic distribution}
$$\chi \coloneqq TS \cap \xi \subset TS$$
in the following manner:
\begin{itemize}
\item each component of the locus where $\chi$ is two-dimensional is endowed with a sign, which is positive if and only if the orientations of $\xi$ and $TS$ coincide;
\item each one-dimensional leaf of $\chi$ is oriented so that, along the leaf, we have an orientation preserving identification
 \begin{equation}
 \label{eq:orientation} (TS/\chi) \oplus (\xi /\chi) \oplus \chi \cong TY = \langle\OP{Reeb}\rangle \oplus \xi,
 \end{equation}
 where $TS/\chi$ and $\xi/\chi$ are oriented according to the rules that that $\chi \oplus (TS/\chi) = TS$ and $\xi /\chi \oplus \chi=\xi$ induce the given orientations on $TS$ and $\xi$, respectively. 
\end{itemize}
Recall the standard fact that the characteristic distribution of a hyperbolic (resp. elliptic) complex tangency of either sign has index $-1$ (resp. +1), where this index refers to the index of the non-degenerate singularity at the point $p$ of any locally defined non-degenerate vector field that is tangent to the distribution. Hence, the following holds: The characteristic distribution in the hyperbolic case corresponds to a saddle with two unstable and two stable manifolds, while the characteristic distribution in the elliptic case corresponds to a focus point with two-dimensional stable or unstable manifold. In order to determine which directions are stable and which are unstable the following lemma is useful:
\begin{lma}
 \label{lma:computingorientations}
 Let $p \in S$ be a point at which $\chi_p =\xi_p \cap T_pS$ is one-dimensional, where $S^2 \subset (Y^3,\alpha)$ is an oriented surface in a contact 3-manifold equipped with a contact form.
 \begin{enumerate}
 \item Assume that $N \in \xi_p$ is an oriented outward normal to $T_pS \neq \xi_p$, and that $R \in T_pS$ satisfies $\alpha(R) > 0$. The vector $V \in \chi_p$ is positively oriented if and only if $\langle V,R\rangle=T_pS$ is an oriented basis.
 \item If $N \in \xi_p$ is the oriented outwards normal to $S$ with the property that complex multiplication gives $JN \in T_pS$, then $V=JN \in \chi_p$ agrees with the orientation of $\chi_p$.
\item A positive (resp. negative) elliptic complex tangency has a two-dimensional unstable (resp. stable) manifold.
 \end{enumerate}
 \end{lma}
 \begin{proof}
 (1): Note that $\langle N,V,R \rangle=T_pY$ is a positively oriented basis by assumption. Hence, so is $\langle R,N,V\rangle=T_pY$ and, thus, $\langle N,V \rangle=\xi_p$ is a positively oriented basis. The statement now follows from Formula \eqref{eq:orientation}.
 
 (2): Again, note that $\langle R, N,V \rangle =T_pY$ is a positively oriented basis by assumption, where $R\in T_pS$ is any vector such that $\alpha(R)>0$. Since $N$ is a positive outward normal and $\langle R,N,V \rangle=\langle N,V,R \rangle$ induce the same orientations, it follows that $\langle V,R \rangle =T_pS$ is a positively oriented basis. The result then follows from Part (1).

 (3): We prove the statement in the positive case; the negative case is analogous. Let $S$ have a positive elliptic complex tangency at $q \in S$ or, equivalently, $d\alpha|_{TS}$ is a positive area form at $q$ and $T_qS=\xi_q$. For two locally defined coordinate vector fields $X,Y$ near $q \in S$ that restrict to a positively oriented basis $\langle X_q, Y_q\rangle T_qS$ we thus compute
 $$ 0<d\alpha(X_q,Y_q)=X_q(\alpha(Y))-Y_q(\alpha(X)) $$
 which implies that at least one of the terms $X_q(\alpha(Y))$ or $-Y_q(\alpha(X))$ has to be positive. Since we can assume that the integral curve of $X$ as well as $-Y$ passing through $q$ is tangent to $\chi$ near $q$, we conclude that at least one of these integral curves agree with the orientation of $\chi$. By symmetry, it follows that this is true for all integral curves emanating from $q$.
 \end{proof}

 \subsection{Unstable and stable manifolds of $h_\pm \in C_\epsilon$.}
 We are now ready to investigate the behavior of the stable and unstable manifolds of the hyperbolic complex tangencies $h_\pm$ contained inside the model cylinder $C_\epsilon$. 
 
\begin{figure}[htp]
\vspace{3mm}
\labellist
\pinlabel $\theta$ at 106 116
\pinlabel $x$ at 221 49
  	\pinlabel $\color{red}h_-$ at 117 75
  \pinlabel $\color{red}h_+$ at 117 36
  \pinlabel $\color{red}h_-$ at 117 -1
	\endlabellist
 \includegraphics{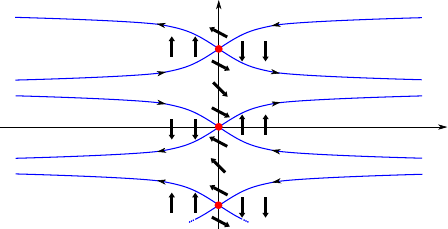}
 \caption{The characteristic distribution $\xi \cap TC_\epsilon$ with the stable and unstable manifolds of $h_\pm$ shown in blue. Here $\theta$ is an $S^1$-valued coordinate on the cylinder, which can be taken to coincide with a translation of the ambient coordinate $\pm y$ near $h_\pm$. For the depicted orientations of the characteristic line field, we have oriented $C_{\epsilon}$ as the boundary of the solid cylinder, i.e.~making $\langle \partial_x,\partial_\theta \rangle$ an oriented basis.}
 \label{fig:sigma1}
\end{figure}

\begin{lma}
\label{lma:tangencies}
The characteristic distribution $TC_\epsilon \cap \xi$ has only two singular points
$$\left\{ (x,y,z)=h_\pm \coloneqq \left(0,0,\pm \sqrt{\epsilon\rho(0)}\right)\right\} \subset \tilde{\Sigma}_\epsilon$$
which are hyperbolic points of opposite orientation signs. When $C_\epsilon$ is oriented as the boundary of the solid cylinder $\{\epsilon^{-1}(y^2+z^2) \le \rho (x)\}$ it follows that $h_+$ (resp. $h_-$) is a positive (resp. negative) hyperbolic complex tangency. Furthermore:
\begin{itemize}
\item The stable as well as unstable manifolds of $h_\pm$ for the characteristic distribution each consists of a pair of integral curves that intersect the boundary of the cylinder $C_\epsilon$ transversely in two points that are contained in different components of $\partial C_\epsilon$;
\item The characteristic distribution $\chi$ along the simple closed curve $\sigma = \{x=0\} \cap C_\epsilon \ni h_\pm$ is tangent to the latter precisely at the singular points $h_\pm$. Moreover, given some preferred normal vector field of $\sigma \subset C_\epsilon$, the oriented line field $\chi$ agrees with the normal orientation along one of the two arcs in $\sigma \setminus \{h_+,h_-\}$, and disagrees along the other one; and 
\item Near the boundary of $\Sigma^U_\epsilon \cap \Sigma^U$ (which coincides with the boundary of the cylinder $C_\epsilon$) the characteristic distribution coincides with $\R\partial_y$ in the above coordinates.
\end{itemize}
 See Figure \ref{fig:sigma1} for a schematic depiction.
\end{lma}
\begin{proof}
The statements can be checked explicitly by investigating the characteristic distribution on $C_\epsilon$; see Figure \ref{fig:sigma1} for a schematic picture of the characteristic distribution and the stable/unstable manifolds.

In addition, note that $(x,y,z)\mapsto(-x,-y,z)$ is a contact-form preserving isomorphism that fixes the cylinder $C_\epsilon$ set-wise, preserves its orientation, while it interchanges the two boundary components, and has the circle $\sigma$ as its fixed-point locus. It immediately follows that this involution preserves the orientation of the characteristic distribution as well. Thus, this involution necessarily preserves the unstable and stable submanifolds.

One immediately concludes that the pairs of integral curves of the stable (resp. unstable) submanifolds of $h_\pm$ intersect the boundary of $C_\epsilon$ in different components.

 The second bullet point is also a consequence of this symmetry.

The last bullet point is immediate, since $\Sigma^U$ is contained in $\{x=\pm 1\}$ in the coordinates for the standard model.
\end{proof}

By the previous lemma, when $\Sigma$ is totally real, the surface $\Sigma_\epsilon$ produced from $\Sigma$ by (possibly several) ambient 0-surgeries along Legendrian arcs with boundary on $\Sigma$ has only hyperbolic complex tangencies. More precisely, there is one pair of hyperbolic points for each attached handle.

\begin{ex}
The Clifford torus
$$T_{\OP{Cl}} \coloneqq S^1_{\frac{1}{\sqrt{2}}} \times S^1_{\frac{1}{\sqrt{2}}} \subset S^3$$
is totally real inside $S^3$ (it is even Lagrangian). The intersection $S^3 \cap \Re{\CC^2}$ is the standard Legendrian unknot, and it intersects the Clifford torus transversely in precisely the four points $\{(\pm1,0),(0,\pm 1)\}$. In other words, this gives four different Legendrian arcs that we can use when performing the surgery construction. An application of the Reeb flow (i.e.\ multiplication by $e^{i\theta}$) produces two $S^1$-families of such arcs; each $S^1$-family consists of pairwise disjoint arcs. Denote by
$$ T_{\OP{Cl}}^g \subset (S^3,\xi).$$ the genus-$g$ surface with precisely $2(g-1)$ hyperbolic points obtained by performing an ambient surgery on $g-1$ arcs of the type described above.
\end{ex}

The remaining part of the section will consist of analyzing holomorphic discs of Maslov index two inside $D^4$ with boundary on the genus-$g$ surface $T^g_{\OP{Cl}} \subset S^3=\partial D^4$. The goal is to establish that $T^g_{\OP{Cl}}$ is rationally convex by showing that these holomorphic discs constitute a holomorphic filling of the surface that satisfies certain additional properties.

\subsection{Holomorphic fillings of a surgery on the Clifford torus}

\label{sec:clifford}

Recall that $T_{\OP{Cl}} = S^1_{1/\sqrt{2}} \times S^1_{1/\sqrt{2}} \subset S^3$ admits two fillings $\mathcal{T}_1$ and $\mathcal{T}_2$ by holomorphic discs of Maslov index two that are contained inside lines; these fillings are the two solid tori
$$\mathcal{T}_1=D^2_{1/\sqrt{2}} \times S^1_{1/\sqrt{2}} \:\:\text{and}\:\: \mathcal{T}_2=S^1_{1/\sqrt{2}} \times D^2_{1/\sqrt{2}}$$
that carry the holomorphic disc foliations
$$ D^2_{1/\sqrt{2}} \times \{e^{i\theta}/\sqrt{2}\}\:\: \text{and} \:\: \{e^{i\theta}/\sqrt{2}\} \times D^2_{1/\sqrt{2}}.$$
For $\epsilon>0$ sufficiently small in the above ambient surgery construction, we may assume that plenty of holomorphic discs in the above two families $\mathcal{T}_i$ have boundaries that are contained inside the subset $T_{\OP{Cl}} \cap T_{\OP{Cl}}^g \subset T_{\OP{Cl}}$ of the Clifford torus that is left undeformed by the surgery. Our goal is to show that these holomorphic Maslov-two disc live in moduli spaces that provide holomorphic fillings of the positive genus surface $T_{\OP{Cl}}^g$ produced by the surgery. Moreover, the fillings produced on the surface of genus $g$ produced by the surgery will be seen to be given by a standard handle-body bounding that surface.

Denote by $\mathcal{T}'_i \subset \mathcal{T}_i$ the union of the discs whose entire boundary is contained inside $T_{\OP{Cl}} \cap T_{\OP{Cl}}^g \subset T_{\OP{Cl}}^g$. For a carefully created surgery in a sufficiently small neighborhood, we may assume that $\mathcal{T}'_i$ is a disjoint union of a number $2(g-1)$ of embedded solid cylinders with boundary on $T_{\OP{Cl}}^g$ of the form
$$\mathcal{T}'_1=D^2_{1/\sqrt{2}} \times A_1 \:\:\text{and}\:\: \mathcal{T}'_2=A_2 \times D^2_{1/\sqrt{2}}$$
where $A_i \subset S^1_{1/\sqrt{2}}$ consists of $2(g-1)$ number of open intervals. When $g>1$ the above surgery may be assumed to produce a surface $T_{\OP{Cl}}^g$ for which
 $$\Sigma_i \coloneqq \overline{T_{\OP{Cl}}^g \setminus \mathcal{T}'_i}$$
 consists of a disjoint union of $g-1$ number of embeddings of the connected compact surface of genus $0$ with four boundary components (i.e.~a sphere with four open balls removed), and for which each connected component is totally real away from precisely two generic hyperbolic complex tangencies. Denote by $\Sigma_i^j$, $j=1,\ldots,g-1$ an enumeration of the connected components of $\Sigma_i$.

\begin{figure}[htp]

\labellist
	\pinlabel $\mathcal{T}^+_i$ at 40 45
	\pinlabel $\mathcal{T}^-_i$ at 40 80
  	\pinlabel $h_-$ at 155 83
	\pinlabel $\mathcal{T}^+_i$ at 275 45
	\pinlabel $\mathcal{T}^-_i$ at 275 80
  \pinlabel $h_+$ at 155 40
	\endlabellist
 \includegraphics{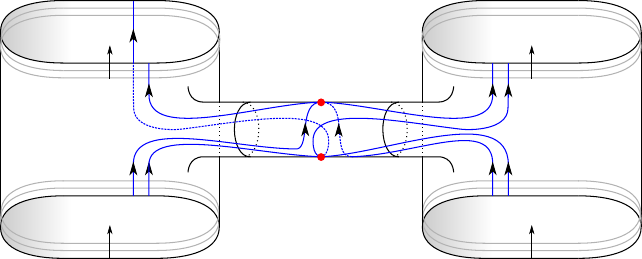}
 \caption{The surface $\Sigma_i^j$ is a sphere with four discs removed. The filling produced by Bedford--Klingenberg yields the handle-body that bounds the surface in the picture, which contains the parts of the disc families $\mathcal{T}^\pm_i \subset \mathcal{T}_i$.}
 \label{fig:sigma2}
\end{figure}
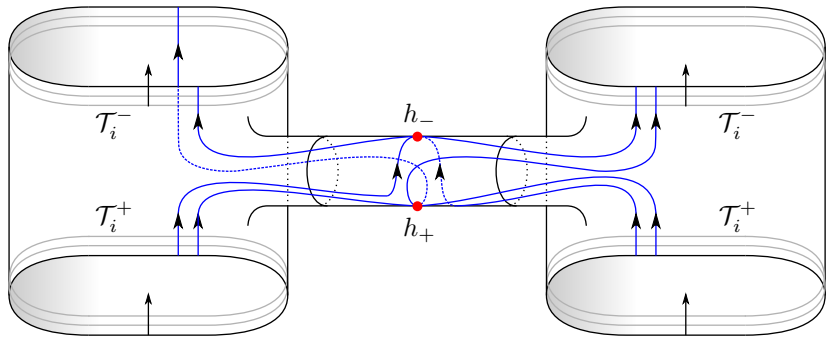

\begin{lma}
 \label{lma:surgery}
 After a generic $C^2$-perturbation of $\Sigma_i$ supported away from the boundary, $\Sigma_i$ admits a holomorphic filling that coincides with $\mathcal{T}_i$ near $\partial \Sigma_i$. Furthermore, the filling of each component $\Sigma_i^j$ is diffeomorphic to a three-ball with four tubes attached (given by the discs in $\mathcal{T}_i)$, where precisely two of the four components of $\partial\Sigma_i^j$ have boundary orientation that coincides with the boundary orientation of the corresponding holomorphic disc in $\mathcal{T}_i$. See Figure \ref{fig:sigma2}.
\end{lma}
\begin{proof}
 Recall the famous result by Bedford--Klingenberg \cite{BedfordKlingenberg} and Kruzhilin \cite{Kruzhilin1992} , by which any generic sphere inside a the boundary of a pseudoconvex four-dimensional domain with only good hyperbolic tangencies admits a filling that is diffeomorphic to a ball. We refer to \cite{BedfordKlingenberg} for the definition of a good hyperbolic tangency, but note that any hyperbolic tangency can be $C^2$-perturbed in an arbitrarily small neighbourhood in order to make it good.

 The main mechanism of the proof of the above existence result for analytic fillings is that the unique one-parameter so-called Bishop families of embedded holomorphic discs that are born near each elliptic complex tangency extend to an analytic filling of the entire sphere. This is done by carefully analyzing the possible breakings at the hyperbolic points, together with delicate use of positivity of intersection for curves with boundary in order to deduce that the discs stay embedded and disjoint. In our setting each $\Sigma_i^j$ is a genus-0 surface with boundary without any elliptic tangencies. However, since the boundary components of $\Sigma_i^j \subset S^3$ all are boundaries of embedded holomorphic discs in $D^4$, that create partial fillings near the boundaries, it is possible to reduce the construction of the filling to the case of the sphere. There are two natural ways to proceed here.

 One alternative, is to complete each boundary component $\partial \Sigma_i^j$ to an embedded sphere by adding discs with exactly one elliptic complex tangency at each of its four boundary components; e.g.~take a perturbation of a radial projection of the holomorphic discs in $\mathcal{T}_i \subset D^4$ to $S^3$. The aforementioned results about fillings by holomorphic discs in four-dimensional symplectic manifolds implies that the spheres produced are fillable by holomorphic discs. A standard argument involving positivity of intersection then shows that this disc family necessarily coincides with the discs from the family $\mathcal{T}_i'$ near the boundary of $\Sigma_i$.

Alternatively, one could simply start with the disc families near $\partial \Sigma_i$ themselves and run the argument of Bedford--Klingenberg to produce a filling, without passing to the auxiliary sphere. Here it is important to use the fact that each of the four boundary components of the surface has a partial filling by discs whose boundaries form collars for the boundary of the surface. \end{proof}

We immediately conclude the following:
\begin{cor}
\label{cor:filling}
After a generic $C^2$-small perturbation of $T^g_{\OP{Cl}}$ support away from the discs in the family $\mathcal{T}_i'$, the latter family extends to a holomorphic filling by discs of $T^g_{\OP{Cl}}$ that is homeomorphic to a handle-body of genus $g$.
\end{cor}

\subsection{Graph structure of the moduli space of discs}

In order to prove rational convexity we need to further analyze the fillings produced by Corollary \ref{cor:filling}. Recall that these fillings are given by the solutions in certain moduli spaces of holomorphic discs with boundary on the surface $ \Sigma_\epsilon$, which exist by the argument by Bedford--Klingenberg \cite{BedfordKlingenberg} and Kruzhilin \cite{Kruzhilin1992}. 

In \cite{BedfordKlingenberg} and \cite{Eliashberg:Filling}, the structure of the moduli space of embedded holomorphic discs of Maslov index two and boundary in an oriented surface was analyzed, and shown to be a directed graph with only vertices of valency one and three for a generic surface. The graph is oriented in the following manner: 
\begin{itemize}
\item The edges parametrize one-parameter families of smoothly embedded Maslov-two discs with boundary on the surface, where the boundary of the discs are disjoint from the complex tangencies, and transverse to the characteristic foliation. Moreover, the edges are oriented so that $dt\wedge d\theta$ agrees with the chosen orientation of the surface under the boundary evaluation map, where $t$ is a coordinate on the oriented edge, and $\theta$ is a coordinate on the boundary of the disc with the orientation induced by the complex structure on the disc.
\end{itemize}
 
This convention makes sense since the embedded discs in the moduli spaces under consideration are embedded and have boundaries that provide local foliations of the surface. 

The orientation of the moduli spaces has the following properties with our conventions:
\begin{lma}
 \label{lma:graphorientations}
Let $S \subset (\partial D,\alpha)$ be a surface with only elliptic and good hyperbolic complex tangencies contained inside the boundary of a smooth pseudoconvex domain $D$. The graph corresponding to the compactified moduli space of a filling of the surface by embedded holomorphic Maslov-2 discs contained inside $D$ satisfies the following properties:
\begin{enumerate}
\item An embedded holomorphic disc $(D^2,\partial D^2) \to (D,S)$ whose boundary is disjoint from the complex tangencies of $S$ is transverse to the characteristic foliation of $S \subset \partial D$ along its boundary. Moreover, $\alpha(\partial_\theta)>0$ is satisfied for the complex boundary orientation, and $\langle V,\partial_\theta\rangle$ is an oriented basis of $TS$ where $V$ denotes the oriented characteristic vector field.
 
\item The one-valent vertices correspond to elliptic complex tangencies and the edge is outgoing if and only if the elliptic point is positive (for $T^g_{\OP{Cl}}$, however, there are no such tangencies); 
\item The three-valent vertices are in bijection with the hyperbolic complex tangencies, where the vertex is a nodal configuration that corresponds to two holomorphic discs that are smooth embeddings away from the node. Moreover, the hyperbolic point is positive (resp. negative) if and only if there are two incoming (resp. outgoing) and one outgoing (resp. incoming) edges.

 As a consequence, the normals in $S$ to the boundary of the holomorphic disc given by, on one hand, the characteristic vector field $V$ and, on the other hand, the infinitessimal variations $\partial_t$ of the oriented moduli space, are both outwards pointing.
\end{enumerate}
\end{lma}
\begin{proof}
 For the property that the moduli spaces of embedded holomorphic Maslov-2 discs has a compactification where one-valent vertices correspond to Bishop families at elliptic complex tangencies, and three-valient vertices correspond to nodal configuration hyperbolic complex tangencies, we refer to \cite{BedfordKlingenberg} or \cite{Eliashberg:Filling}.

 (1): The property that embedded holomorphic discs with boundary disjoint from the complex tangencies on $S \subset \partial D$ are transverse to the characteristic foliation is a standard property of holomorphic discs with boundary on totally real submanifolds contained inside a pseudoconvex domain; see e.g.~\cite{Zehmisch2015}. The property that the positively oriented characteristic vector field is a positively oriented normal then follows directly from Part (1) of Lemma \ref{lma:computingorientations}.
 
 (2): This follows directly from Part (3) of Lemma \ref{lma:computingorientations}.
 
 (3): This statement can be verified in a suitable standard model. For instance, one can consider a smooth $S$ sphere inside the Levi-flat hypersurface
 $$\C_{x_1+iy_1} \times \R_{x_2} \subset \C_{x_1+iy_1} \times \C_{x_2+iy_2},$$
 where the orientation of this hypersurface is given by $\langle \partial_{x_1},\partial_{y_1},\partial_{x_2}\rangle$, and where the sphere has been oriented as the boundary of an open 3-ball. One can construct such a sphere with precisely the following complex tangencies: Two positive elliptic complex tangencies at $x_2=2$, two negative complex elliptic tangencies at $x_2=-2$, one positive hyperbolic complex tangency at $x_2=1$, and one negative hyperbolic complex tangency at $x_2=-1$. The sphere $S$ has an obvious filling where the discs are given by the compactifications of the bounded domains in the slices $\{x_2=\OP{cst}\} \setminus S$.

 In this example, the two Bishop families born at the positive elliptic points join at the hyperbolic point at $x_2=1$, then split at the hyperbolic point at $x_2=-1$, and finally die at the two negative elliptic point at $x_2=-2$. The positive (resp. negative) hyperbolic complex tangency at $x_2=1$ (resp. $x_2=-1$) has two incoming (resp. outgoing) edges and one utgoing (resp. incoming) edge.

 Since this is a model for the general filling, the same property holds for a generic 3-valent node of the moduli space; see Figure \ref{fig:split}.
 \end{proof}

\begin{figure}[htp]
 \label{fig:split}
 \includegraphics{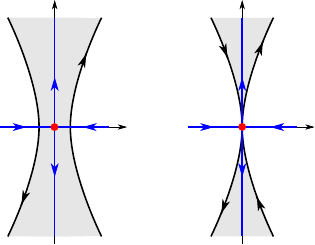}
 \caption{The model a single incoming edge at a hyperbolic point, with two outgoing edges. Note that the hyperbolic point is negative here. (Using the terminology from \cite{BedfordKlingenberg}, the discs approach the hyperbolic point ``from the outside.'')}
 \label{fig:split}
\end{figure}

 Note that the discs at the incoming edge(s) approach the hyperbolic point \emph{from the outside} (resp. \emph{inside}) when the hyperbolic point is negative (resp. positive), according to the notation of \cite{BedfordKlingenberg}. See Figure \ref{fig:split}.

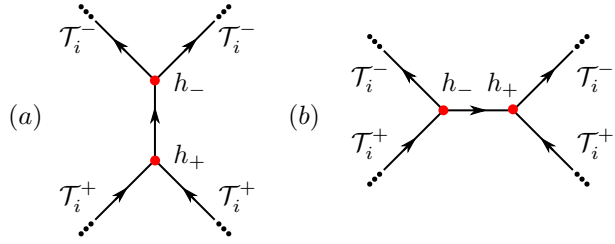
\begin{figure}[htp]
 	\labellist
	\pinlabel $\mathcal{T}^+_i$ at 60 15
	\pinlabel $\mathcal{T}^-_i$ at 60 75
  	\pinlabel $\mathcal{T}^-_i$ at 0 75
 	\pinlabel $h_-$ at 42 57
 	\pinlabel $h_+$ at 42 30
 	\pinlabel $h_-$ at 143 57
 	\pinlabel $h_+$ at 160 57
	\pinlabel $\mathcal{T}^+_i$ at 110 35
  	\pinlabel $\mathcal{T}^-_i$ at 110 60
  \pinlabel $\mathcal{T}^-_i$ at 195 60
 	\pinlabel $\mathcal{T}^+_i$ at 195 35
	\pinlabel $\mathcal{T}^+_i$ at 0 15
		\pinlabel $(a)$ at -20 45
				\pinlabel $(b)$ at 85 45
	\endlabellist
 \includegraphics{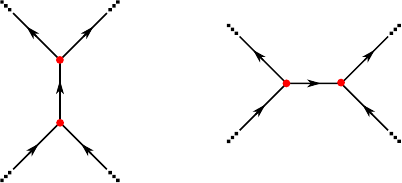}
 \caption{The a priori possibilities of the graphs that corresponds to the moduli space of the discs produced by \cite{BedfordKlingenberg,Kruzhilin1992} and that fill a component $\Sigma_i^j \subset \Sigma_i$. Configuration (a) is shown on the left, where the disc family $\mathcal{T}^{j,+}_i$ form the two incoming edges at $h_+$, while the disc family $\mathcal{T}^{j,-}_i$ form the two outgoing edges at $h_-$. Configuration (b) is shown on the right, where one of the two disc families in $\mathcal{T}^{j,+}_i$ forms the unique incoming edge at $h_-$, while one of the two disc families in $\mathcal{T}^{j,-}_i$ forms the unique outgoing edge at $h_+$.}
 \label{fig:moduli}
\end{figure}

\begin{lma}
 \label{lma:graph}
The graphs corresponding to the moduli spaces that fill a component $\Sigma_i^j \subset \Sigma_i$ can a priori be either of Configuration (a) or (b) shown in Figure \ref{fig:moduli}, where (a) consists of two incoming (resp. outgoing) families from the discs in $\mathcal{T}_i'$ at the positive (resp. negative) hyperbolic point.
\end{lma}
\begin{proof}
 
 The filling of $\Sigma^j_i$ given by Lemma \ref{lma:surgery} consists of discs from $\mathcal{T}_i'$ whose boundaries foliate the collars of the four boundary components $\partial \Sigma^j_i$ by construction. Lemma \ref{lma:graphorientations} implies that two of these components correspond to edges oriented into the filling of $\Sigma^j_i$, which we denote by $\mathcal{T}_i^{j,+}$, while two are oriented out from $\Sigma^j_i$, which we denoted by $\mathcal{T}_i^{j,-}$.

It is now simply a matter of enumerating all directed graphics with two three-valent vertices, two incoming leaves, and two outgoing leaves.
\end{proof}

\begin{figure}[htp]
 	\labellist
 	\pinlabel $\mathcal{T}_i'$ at 28 61
  \pinlabel $\mathcal{T}_i'$ at 28 28
 	\pinlabel $\mathcal{T}_i'$ at 67 61
  \pinlabel $\mathcal{T}_i'$ at 67 28
 	\pinlabel $\mathcal{T}_i'$ at 105 61
  \pinlabel $\mathcal{T}_i'$ at 105 28
 	\pinlabel $\mathcal{T}_i'$ at 142 61
  \pinlabel $\mathcal{T}_i'$ at 142 28
 	\pinlabel $\mathcal{T}_i'$ at 191 61
  \pinlabel $\mathcal{T}_i'$ at 191 28
	\endlabellist
  \includegraphics{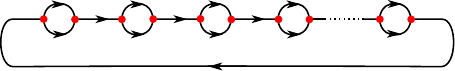}

 \caption{The filling of $T^g_{\OP{Cl}}$ obtained by gluing together the partial fillings shown on the left in Figure \ref{fig:moduli} for our choices of Legendrian arcs for the ambient surgery that produces $T_{\OP{Cl}}$.}
 \label{fig:handlebody}

\end{figure}

In order to prove rationally convexity of $T^g_{\OP{Cl}}$ we need to sharpen the result of Corollary \ref{cor:filling} to the following structural result by excluding Configuration (b) from appearing; see Figure \ref{fig:moduli}.
\begin{prp}
\label{prp:filling}
The holomorphic filling by disks of $T^g_{\OP{Cl}}$ produced by Corollary \ref{cor:filling} is a genus-$g$ handle body that contains the disc family $\mathcal{T}_i'$ as a sub-filling, and for which
\begin{itemize}
\item the vertices in the moduli space are all three-valent and in bijective correspondence with the $2(g-1)$ number of hyperbolic complex tangencies;
\item each pair of hyperbolic complex tangencies living in the four-punctured sphere $\Sigma_i^j \subset T^g_{\OP{Cl}}$ is contained in a sub-graph of the type Configuration (a), as shown in Figure \ref{fig:moduli}; in particular, any oriented edge in the moduli space can be extended to an oriented cycle in the graph.
\end{itemize}
The moduli space is the graph shown in Figure \ref{fig:handlebody}; note that the pair of hyperbolic points that are contained inside the same handle are connected by a pair of oriented edges that share the same start an endpoints.
\end{prp}
\begin{proof}
 Lemma \ref{lma:graph} implies that the filling of $\Sigma_i^j$ corresponds to one of the two graphs shown in Figure \ref{fig:split}. The main goal is to exclude Configuration (b) on the right in the figure. This will be done by a careful consideration of the characteristic foliation on $\Sigma_i^j$.

After having ruled out Configuration (b), it is not hard to see that any edge can be joined to a cycle, as claimed in the second bullet point. Namely, we first contract the edges in the graph that connect the pairs of hyperbolic complex tangencies from each $\Sigma_i^j$. The new graph has now only four-valent vertices, each with two incoming and two outgoing edges. Since this graph can be decomposed into a union of cyclic graphs merged at certain vertices, it is obvious that the new graph satisfies the property that any edge can be completed to a loop. It is not hard to see that the same property also holds for the original graph in the case when each contracted edge is of the type shown in Configuration (a). (This is, however, not necessarily the case if there are edges as in Coniguration (b).)

We then proceed to exclude Configuration (b). Recall that each component $\Sigma_i^j$, $i=0,1$, $j=1,\ldots,2g-1$ contains the standard model of a cylinder $C_\epsilon$ from Section \ref{sec:fillable}. Denote this cylinder by $C_i^j \subset \Sigma_i^j$ and the pair of hyperbolic points $(h_i^j)^\pm \in C_i^j$. 

 Lemma \ref{lma:tangencies} shows that the two components of the unstable (resp. stable) manifold of the critical points $(h_i^j)^\pm$ of the characteristic foliation intersect the boundary $\partial C_i^j$ in two different components, where each component intersects the boundary transversely in a single point. Using the fact that the characteristic distribution of $\mathcal{T}_{\OP{Cl}} \subset S^3$ is tangent to the vector-field generated by the standard action of $S^1$ by scalar multiplication on $\C^2$, which is everywhere transverse to the standard disc family $\mathcal{T}_i$, we conclude that each component of the unstable (resp. stable) manifolds of each of complex tangency $(h_i^j)^\pm$ intersects a unique component of the boundary of $\Sigma_i^j$ transversely in a single point. (Recall that $\partial \Sigma_i^j$ consists of four bounday components, along two of which the characteristic field is inwards and outwards pointing, respectively.) 

This means that the discs in the family $\mathcal{T}_i^{j,+} \subset \mathcal{T}'_i$, whose boundaries foliate a neighborhood of the collars of the boundary of $\Sigma_i^j$, such that the orientation of the moduli space of $\Sigma_i^j$ is inwards-pointing, cannot be a part of a unique incoming edge of the negative hyperbolic point $(h_i^j)^-$. Or, in other words, this disc family cannot split at a hyperbolic point (see Lemma \ref{lma:graph}). Namely, as can be seen by studying the behavior of such a disc family at the moment close to splitting (see Figure \ref{fig:split}), the disc must necessarily intersect \emph{both} components of the stable manifold of $(h_i^j)^-$ at that moment. Since the holomorphic discs remain transverse to the characteristic distribution as long as they do not intersect the complex tangencies (see Lemma \ref{lma:graph}), this leads to a contradiction with the previously established fact that the discs near the boundary of $\Sigma_i^j$ only intersect one of the components of the stable manifold. 

 This excludes the possibility of Configuration (b), and we are thus left with Configuration (a) as sought.
\end{proof}

\subsection{Holomorphic fillings of a surgery on a knotted rationally convex torus}
\label{sec:knotted}

In this section we provide two approaches for constructing rational convex high genus surfaces in $S^3$ whose smooth isotopy class inside the sphere differs from the $T^g_{\OP{Cl}}$ constructed in Section \ref{sec:clifford}. Note that the latter surface is unknotted in the sense that it creates a Heegaard-splitting of $S^3 \setminus T^g_{\OP{Cl}}$ into two handle-bodies.

 Take any transverse knot $K \subset S^3 \subset \C^2$. Recall the standard fact that such a knot can be constructed in any smooth isotopy class. The standard neighborhood theorem of transverse knots implies that it has a solid torus neighborhood which is contactomorphic to
$$(S^1_\phi \times D^2_\epsilon, d\phi+r^2/2\,d\theta),$$
under which $K$ is identified with the core $S^1 \times \{0\}$; here $(r,\theta)$ denotes polar coordinates on the small disc $D^2_\epsilon$. Note that this identification does not identify contact forms, but that the identity $\alpha_{st}=e^g(d\phi+r^2/2\,d\theta)$ holds for some smooth function $g$ defined on the above neighborhood.

For any smooth function
$$f \colon S^1 \times S^1 \to (0,\epsilon)$$
we will consider the totally real torus defined by
$$ T_f=\{ r=f(\phi,\theta) \} \subset S^1_\phi \times D^2_\epsilon.$$
Note that $T_f$ is pre-Lagrangian when $\partial_\phi f \equiv 0$ is satisfied.

Any torus $T_f$ has at least one complement in $S^3$ that is a solid torus, namely 
$$S_f=\{ r \le f(\phi,\theta) \} \subset S^1_\phi \times D^2_\epsilon$$
whose core is the original transverse knot $K$. Moreover, the solid torus $S_f$ can be seen to be foliated by symplectic discs
$$D_{f}(\alpha) \coloneqq S_f \cap \{\phi=\alpha\} \subset S^1_\phi \times D^2_\epsilon , \:\:\alpha \in S^1,$$
having boundary on $T_f$ and being of Maslov index equal to two. In the case when the original transverse knot $K$ is knotted, the other remaining component is not a solid torus. In this case, $T_f \subset S^3$ is of course also knotted.

Next we exhibit families of totally real tori of the form $T_f$ that admit fillings by holomorphic discs with boundary in the class of the meridian, and which bound smooth solid tori in $S^3$ whose core is the transverse knot $K$. Note that the existence of the filling implies that these tori are rationally convex; see \cite[Lemma in 1.(c)]{Duval:Riemann}.
\begin{lma}
 \label{lma:foliationbyfillable}
 A neighbourhood of $K$ is foliated by totally real tori of the form $T_{f_t}$ that admit embedded pseudoholomorphic disc fillings, when $f_t$ have been suitable chosen. We can assume that $\partial D_{f_t}(0) \subset T_{f_t}$ is a boundary of one of the holomorphic discs in the filling. In addition, the functions $f_t$ can be chosen to satisfy e.g. one of the following natural conditions;
 \begin{itemize}
 \item $f_t=t$ away from some arbitrarily small neighbourhood of $\{\phi=0\}$; or
 \item $\partial_\phi f_t=0$ holds everywhere, i.e.~making the tori $T_{f_t}$ pre-Lagrangian.
 \end{itemize}
\end{lma}
\begin{proof}
 Consider first the disk $\{\phi=0\} \times D^2_\epsilon$ whose characteristic distribution has an isolated elliptic singulatity at $(0,0)$. There is a Bishop family born at this elliptic complex tangency, which near the elliptic point creates a foliation by boundaries of holomorphic discs. We let $g_t(\theta)$, $t \ge 0$, be functions such that the graph $\{r=g_t(\theta)\} \subset D^2_\epsilon$ defines this foliation, normalised by $g_t(0)=t$, where $g_0 \equiv 0$ is the unique singular leaf corresponding to the elliptic complex tangency. Denote by $D_t \subset D^4$ the holomorphic disc with boundary $\{r=g_t(\theta)\}$.

 Let $f_t \colon S^1 \times S^1 \to (0,\epsilon)$ be any smooth function that satisfies $f_t(0,\theta)=g_t(\theta)$. Using Zehmisch result from \cite{Zehmisch2015} we conclude that the totally real torus
 $$T_{f_t} \coloneqq \{r=f_t(\phi,\theta)\} \subset \left(S^1_\phi \times D^2_\epsilon,d\phi+r^2/2\theta\right), \:\: t >0,$$
 which contains the boundary of the embedded holomorphic disc $D_t$ by construction, admits a filling by discs with boundary in the class of the meridian $[\{\phi=\OP{cst}\}] \in H_1(T_{f_t}).$ This filling moreover contains the disc $D_t$ as one of its disc leaves.

 The smooth dependence on the discs on the boundary condition and positivity of intersection implies the claim, after the family of functions $f_t$ have been appropriately chosen so that the tori $T_{f_t}$ creates a smooth foliation of a neighbourhood of $K$.
\end{proof}
\begin{rmk}
 With a bit of more work one can also deduce that the pre-Lagrangian tori $T_\rho=\{r=\rho\}$ for the constant functions $f(\phi) \equiv \rho$ admit fillings by holomorphic discs. However, we do not need this stronger statement for our purposes. We sketch two alternative proofs of this fact:
 \begin{itemize}
 \item The first alternative follows the technique invented by Duval--Gayet in \cite{Duval:Riemann} for constructing a holomorphic disc filling of an arbitrary unknotted rationally convex torus in $S^3$. The construction takes place in the universal cover of a subset $\Omega$ foliated by the holomorphic fillings by the tori provided by Lemma \ref{lma:foliationbyfillable}. In the universal cover, the cylinder $\tilde{T}_\rho \to T_\rho$ that lifts the torus is pinched to provide a very long sphere which can be filled by the technique due to Bedford--Klingenberg. To conclude that the filling also provides discs that have boundary on the original cylinder $\tilde{T}_\rho$ needs a compactness argument that is significantly easier in this case compared to the argument in \cite{Duval:Riemann}, since we have full control of the characteristic foliation of the involved torus, and can thus provide precise area bounds.
 \item The second alternative is to provide a deformation of a holomorphic fillable totally real torus $T_{f_\epsilon}$ provided by Lemma \ref{lma:foliationbyfillable} through totally tori of the form $T_{f},$ until arriving at $T_\rho$, and arguing that the filling persists under this deformation. To that end we also need an a priori control of area, which can be readily obtained in view of the explicit form for the characteristic foliation (especially in the pre-Lagrangian case).
 \end{itemize}
 \end{rmk}

The rest of the construct can be performed as in Section \ref{sec:clifford}, yielding an embedded genus-$g$ surface $T^g \subset S^3$ with precisely $2(g-1)$ number of hyperbolic complex tangencies, which admit holomorphic fillings by discs with the same properties as those in Proposition \ref{prp:filling}. Note that the Legendrian arcs needed for the surgery can be constructed explicitly e.g.~inside the solid torus neighborhood of the transverse knot.

Alternatively, one can also start with the surface given by the Clifford torus as above, but then choose the Legendrians arcs along which the surgery is performed to be knotted.

\subsection{Rational convexity of $T^g$ (Proof of Theorem \ref{thm:ratconvexS3})}
\label{proof:ratconvexS3}

Here we prove the rational convexity of the fillable genus-$g$ surfaces $T^g$ of the type constructed in either of Sections \ref{sec:clifford} and \ref{sec:knotted}.

This argument follows as the proof of rational convexity for fillable tori in \cite[Section 1.b]{Duval:Riemann}. We repeat it here for completeness. The first step is to show that the filling itself is rationally convex.

\begin{lma}
\label{lma:fillingconvex}
The filling $\mathcal{T} \subset D^4$ of $T^g\subset S^3$ as constructed in Subsections \ref{sec:clifford} or \ref{sec:knotted} is rationally convex. In particular, it contains the rational hull of $T^g$, i.e.~$r(T^g) \subset \mathcal{T} \subset D^4$.
\end{lma}
\begin{proof}
 By pushing $T^g$ off of itself along a non-vanishing normal direction inside $S^3$, we can create an open neighborhood of $T^g \subset S^3$ that is foliated by smoothly varying genus-$g$ surfaces, each of which has only isolated good hyperbolic complex tangencies.

 In order to guarantee that the hyperbolic complex tangencies of the push-offs remain good, we utilize the holomorphic $S^1$-action induced by unitary scalar multiplication in $\C^2$ to create the push-offs near the hyperbolic complex tangencies. More precisely, applying this $S^1$-action to the local sheets of $T^g$ that contain the hyperbolic complex tangencies, we get a local foliation by pieces of surfaces that all have good hyperbolic complex tangencies (recall that the $S^1$-action on $S^3$ is transverse to any complex tangency of $T^g$). Extending these push-offs to a push-off of the entire surface $T^g$, we can create a smooth one-parameter foliation by surfaces containing $T^g$ as a leaf, of which each leaf has only good hyperbolic complex tangencies.

It now follows by stability of the solutions, that the fillings persist for the genus-$g$ surfaces in the foliation sufficiently close to the original surface. Positivity of intersection implies that two fillings for two disjoint surfaces in the family are disjoint. In conclusion, these fillings foliate an open neighborhood of $\mathcal{T} \subset D^4$.

We have shown that any point in $D^4 \setminus \mathcal{T}$ that lives in a sufficiently small neighborhood $ U \subset D^4$ of $\mathcal{T}$ lies on a holomorphic disc that is disjoint from the filling $\mathcal{T} \supset T^g$, where this disc has boundary contained in $S^3$. By a standard argument, this Riemann surface can be approximated arbitrarily well by the zero-set of a polynomial that does not vanish on the filling $\mathcal{T}$. After a perturbation, one can moreover construct such a global analytic variety that passes through any fixed chosen point contained in $U \setminus \mathcal{T}$.

It follows from the above that $\mathcal{T}$ is a connected component of the rational hull of $T^g$. From the well-known fact that the rational hull of a connected set again is connected (this is a direct consequence of \v{S}ilov's idempotent theorem), we can finally conclude that $\mathcal{T}$ must contain the entire rational hull, as sought.
\end{proof}

Then we finish by producing analytic varieties that pass through the discs in the filling, and which do not pass through the surface in $S^3$, thus excluding the interior of the discs in the filling from the rational hull.

\begin{lma}
Any point inside the interior $\mathcal{T} \cap B^4$ of the filling of $T^g$ constructed in Subsections \ref{sec:clifford} or \ref{sec:knotted} is not contained inside the rational hull $r(T^g) \subset D^4$ of the boundary of the filling.
\end{lma}
\begin{proof}
The holomorphic filling $\mathcal{T}$ of $T^g$ constructed in the aforementioned subsections, in particular Proposition \ref{prp:filling}, has the following property: Any point $p \in \mathcal{T} \setminus T^g$ contained inside the filling which is disjoint from the boundary $T^g$, lies in a smoothly embedded closed curve that meets every disc in the filling transversely in the interior. It causes no constraint to assume that this curve is real-analytic. Thus, the curve can be extended to a small holomorphic annulus that intersects the filling near the core of the annulus.

Since the filling $\mathcal{T}$ is rationally convex by Lemma \ref{lma:fillingconvex}, we can approximate it by a smooth rationally convex pseudoconvex domain. Then we can run an argument similar to the one in the proof of the previous lemma. First, since the extended domain is a domain of holomorphy, the annulus is cut out by a holomorphic function on the domain. Second, this holomorphic function can be approximated by a rational function on all of $\C^2$. Finally, this means that the annulus can be perturbed to a properly embedded algebraic curve in $\C^2$ that is disjoint from $T^g$. Since the approximation can be assumed to be arbitrarily close, it causes no restriction to assume that this global analytic variety passes throught the point $p$ while remaining disjoint from $T^g$. In other words, $p \notin r(T^g)$, as sought. See \cite[Lemma in 1.(c)]{Duval:Riemann} for a similar argument.
\end{proof}

\section{Symplectic condition for rational convexity}

A famous result by Duval--Sibony \cite[Theorem 3.1]{DuvalSibony} states that any totally real half-dimensional submanifold $\Sigma \subset \C^n$ is rationally convex if and only if it is Lagrangian for a global K\"{a}hler form on $\C^n$ (which in fact can be taken to be equal to the standard linear $\omega_0$ outside of a compact subset). Here we provide a generalisation of this result in the case $n=2$ to surfaces $\Sigma \subset \C^2$ that have a finite number of complex hyperbolic tangencies in addition to the totally real locus. Recall that generic complex tangencies are either hyperbolic or elliptic, and that elliptic complex tangencies are obstructions to rational convexity.

A surface that has a complex tangency can of course not be Lagrangian for a global K\"{a}hler form on $\C^2$. However, the condition that ensures rational convexity is, roughly speaking, that the surface is Lagrangian for a K\"{a}hler form that degenerates precisely at the complex tangencies. To that end, we use Shafikov--Sukhov's adaptation \cite[Lemma 2]{ShafikovSukhov} of Gayet's condition for rational convexity from \cite{Gayet}. Another crucial ingredient that we rely on is the local Stein neighborhood basis of a surfaces with only flat hyperbolic complex tangencies, as constructed by Slapar \cite{Slapar}.

\begin{thm}
\label{thm:rationalconvexity}
Let $\Sigma \subset \C^2$ be a smooth compact surface that is totally real outside of a finite number of flat hyperbolic tangencies $H \subset \Sigma$. Assume that $\Sigma$ is Lagrangian for a global K\"{a}hler form $\omega$ on $\C^2 \setminus H$ that extends to a global smooth closed two-form on $\C^2$, and which in a neighborhood $U_H \supset H$ of the hyperbolic complex tangencies $H$ is given by $\omega=i\del\delbar\rho$ for a function $\rho$ that satisfies
\begin{itemize}
\item $\rho \ge 0$, and $\rho$ vanishes precisely at $H$;
\item $d\rho \neq 0$ in $U_H \setminus H$; and
\item $d^c\rho|_{T\Sigma}$ vanishes on $\Sigma \cap U_H$.
\end{itemize}
Then $\Sigma$ is rationally convex.
\end{thm}

 Before proving the above theorem, we need the following intermediary lemma which creates a deformation $\tilde{\rho}$ of $\rho$ that is pluriharmonic in some neighborhood of $H$. 

\begin{lma}
 \label{lma:phitilde}
 Consider the two-form $\omega$ which satisfies the assumptions of Theorem \ref{thm:rationalconvexity}, and let $U_H \supset H$ again denote the neighborhood of the hyperbolic complex tangencies $H$. For any smaller neighborhood $U_H' \subset U_H$ of $H$, there exists a global plurisubharmonic function $\tilde{\rho} \colon \C^2 \to \R_{\ge0}$ that
 \begin{itemize}
\item is pluriharmonic in some closed neighbourhood $V \subset U'_H$ of $H$;
\item is strictly plurisubharmonic in $\C^2 \setminus V$ and satisfies $i\del\delbar\tilde{\rho}=\omega$ in $\C^2 \setminus U_H$; and
\item satisfies the property that $i\del\delbar \tilde{\rho}|_{T\Sigma}=0$ along $\Sigma$.
\end{itemize}
 \end{lma}
\begin{proof}
Consider the non-negative plurisubharmonic function $\rho \ge0$ defined near $H=\rho^{-1}(0)$ with the properties prescribed by Theorem \ref{thm:rationalconvexity}. In particular, $\omega=i\del\delbar\rho$ holds on $U_H$.

Post-composing the plurisubharmonic function $\rho$ with a suitable convex function $\chi_\epsilon \colon \R_{\ge 0} \to \R_{\ge 0}$ that vanishes in $[0,\epsilon]$, we obtain a plurisubharmonic function $\chi_\epsilon\circ \rho$ that 
\begin{itemize}
\item vanishes in the closed neighborhood $V=\rho^{-1}(-\infty,\epsilon]$ of the hyperbolic points, where we may assume that $ V \subset U'_H$ after taking $\epsilon>0$ sufficiently small;
\item is strictly plurisubharmonic in $ U'_H \setminus V$, while $i\del\delbar(\chi_\epsilon\circ \rho)=\omega$ moreover holds outside of a compact subset of $ U'_H \supset V$;
\item still satisfies the property that $d^c(\chi_\epsilon \circ \rho)$ vanishes near the hyperbolic points when restricted to $T(\Sigma \cap U_H)$.
\end{itemize}
Denote by $\tilde{\omega}$ the globally defined smooth $(1,1)$-form that coincides with $i\del\delbar(\chi_{ \epsilon} \circ \rho)$ in $U_H$ and with $\omega$ in $\C^2 \setminus U_H$. The $i\del\delbar$-lemma then provides the sought global plurisubharmonic function $\tilde{\rho}$ for which $i\del\delbar\tilde{\rho}=\tilde{\omega}$.
\end{proof}

\begin{proof}[Proof of Theorem \ref{thm:rationalconvexity}]
By the Slapar's result \cite[Theorem 2]{Slapar}, after shrinking the neighborhood $U_H \supset H$ of the hyperbolic complex tangencies, we can find a plurisubharmonic function $\psi \colon U \to \R_{\ge 0}$ defined in some neighborhood $U \supset U_H$ of $\Sigma$ that has the property that $\psi^{-1}(0)=\Sigma=d\psi^{-1}(0)$, and where $\psi$ is strictly plurisubharmonic away from the hyperbolic tangencies $H$. We extend $\psi$ to a smooth and compactly supported function defined on all of $\C^2$.

Consider the global plurisubharmonic function $\tilde{\rho} \colon \C^2 \to \R_{\ge0}$ produced by Lemma \ref{lma:phitilde}, and consider the function
$$\phi\coloneqq C\cdot \tilde{\rho}+\psi \colon \C^2 \to \R_{\ge 0}, \:\: C>0.$$
By construction, $\tilde{\rho}$ is pluriharmonic in a neighborhood $V \subset U_H \subset U$ of the hyperbolic tangencies and strictly plurisubharmonic in $\C^2 \setminus V$. If we take $C \gg 0$ sufficiently large, the plurisubharmonic $\phi$ moreover becomes strictly plurisubharmonic on all of $\C^2 \setminus H$.

Rational convexity is then a consequence of \cite[Lemma 2]{ShafikovSukhov} or \cite[Lemme 1]{Gayet}. We proceed to show that the necessary conditions are satisfied, so that he result can be applied.

A standard argument shows that, after perturbing $\phi$ away from $U_H$, while keeping the Lagrangian property of $\Sigma$, we may further assume that $d^c\phi|_{T\Sigma}=dg$ for some smooth $g \colon \Sigma \to \R/2\pi\Z$; see e.g.~\cite{ShafikovSukhov} or \cite{Gayet}. We then consider the function $e^{\phi+ig}$ defined on $\Sigma$. 

First, we claim that $e^{\phi+ig}=e^f$ is satisfied along $\Sigma \cap V$ for some holomorphic function $f \in \mathcal{O}(V)$. Indeed, $\phi=C\tilde{\rho}$ holds along $\Sigma$, since the latter is a critical manifold of the function $\psi$. Since $\tilde{\rho}$ is pluriharmonic in $V$, we can write $C\tilde{\rho}=\Re f$ there, where $f \in \mathcal{O}(V)$. Further, $g=\Im f$ holds up to a constant, since $dg=d^c\phi=d^c C\tilde{\rho}$ along $\Sigma$.

Second, we use Wermer--H\"{o}rmander's result \cite[Lemma 4.3]{WermerHormander} to extend the holomorphic function $e^f$ defined on $V$ to a function $h$ defined in a neighborhood of $\Sigma \cup V$, where further 
\begin{itemize}
\item $h|_\Sigma=e^{\phi+ig}$; and
\item $\delbar h =\mathbf{O}(d(\cdot,\Sigma))^k$ for some arbitrary $k \gg 0$.
\end{itemize}
Since $d^c\phi=dg$ is satisfied along $\Sigma$ by construction, the second bullet point above implies that the functions $|h|$ and $e^\phi$ agree to the first order along $\Sigma$.

Third, after cutting $h$ off by a bump function supported in some small neighborhood of $\Sigma$ the inequality $0<|h| \le e^\phi$ can be assumed to hold globally, with equality precisely along $\Sigma$. Indeed, inside $V$ we have $\log{|h|}=C\cdot\tilde{\rho}$ holds on all of $V$, and since $\psi$ is positive away from $\Sigma$ in $U$, we thus get $|h| \le e^\phi$ in $V \subset U$ with equality precisely on $\Sigma \cap V$. To show the inequality in a sufficiently small neighborhood $U'$ of $\Sigma$ we argue as follows. The function $\phi-\log |h|$ can be assumed to be strictly plurisubharmonic on $U' \setminus H$. Note that, by construction, this plurisubharmonic function vanishes along $\Sigma$ to the first order. After shrinking $U'$ further, a standard argument implies that $|h| \le e^\phi$ holds everywhere in $U'$, with equality precisely along $\Sigma$. Indeed, $\phi-\log |h|$ must have a critical locus along $\Sigma \setminus H$ with a positive definite Hessian in the normal direction. 

This establishes all assumptions needed in order to apply \cite[Lemme 1]{Gayet}, which then shows that $\Sigma$ is rationally convex.
\end{proof}

\subsection{The converse statement: rationally convex are singular Lagrangians}

We then prove a converse to Theorem \ref{thm:rationalconvexity}, by showing that any rationally convex surface with standard flat hyperbolic tangencies can be made Lagrangian for a symplectic form that becomes degenerate precisely at the hyperbolic points.

\begin{thm}
 Assume that $\Sigma \subset \C^2$ is a smooth surface that is rationally convex and totally real except for a finite number of standard flat hyperbolic complex tangencies. There exists a smooth $(1,1)$-form $\omega$ on $\C^2$ which is K\"{a}hler on $\C^2 \setminus H$ for which $\Sigma \setminus H \subset \C^2 \setminus H$ is Lagrangian.
\end{thm}

\begin{proof}
The proof is completely analogous to the case when $\Sigma$ is totally real, which is due to Duval--Sibony \cite[Theorem 3.1]{DuvalSibony}.

The rational convexity of $\Sigma$ implies that there is a smooth function $\phi \colon \C^2 \to \R$ for which $i\del\delbar \phi$ is positive on $\C^2 \setminus \Sigma$, while it vanishes precisely on $\Sigma$. See \cite[Remark 2.2]{DuvalSibony}.

The sought K\"{a}hler form can then be taken to be $\omega=i\del\delbar(\phi+\epsilon \psi)$ where $\epsilon>0$ is sufficiently small, and $\psi$ is the plurisubharmonic function defined in some neighbourhood of $\Sigma$ that was constructed by Slapar in \cite[Theorem 2]{Slapar}. To that end, note that $\psi$ can be assumed to be a global compactly supported smooth function (that, however, only is strictly plurisubharmonic close to $\Sigma$).
 \end{proof}

\section{Lagrangian model of a complex hyperbolic tangency}
\label{sec:symplecticdescription}

The goal of this section is to construct a local model of a hyperbolic complex tangency in $\C^2$ together with a K\"{a}hler form that is degenerate precisely at the tangency, and for which the local model becomes Lagrangian. This construction will then be used as a building block for surfaces that satisfy the assumptions of Theorem \ref{thm:rationalconvexity}, and which thus are rationally convex.

The hyperbolic tangency that we consider is of the particular form
$$ \left\{v=(\Re{u})^2-(\Im u)^2\right\} \subset \C^2_{u,v}$$
in local holomorphic coordinates $(u,v)$ on $\C^2$. We will call such points {\bf standard hyperbolic complex tangencies}. Note that this hyperbolic complex tangency is not ``good'' in the sense of \cite{GaussierSukhov}, but it can be approximated by hyperbolic tangencies that are good. More importantly, this complex hyperbolic tangency is ``flat'' in the sense of Slapar \cite{Slapar}.

\subsection{Local symplectic description of the hyperbolic tangency}
\label{sec:model}
\begin{figure}[htp]
 \vspace{3mm}
 	\labellist
  \pinlabel $y$ at 1 90
  \pinlabel $z$ at 1 195
  \pinlabel $x$ at 223 42
  \pinlabel $x$ at 223 105
	\endlabellist
  \includegraphics{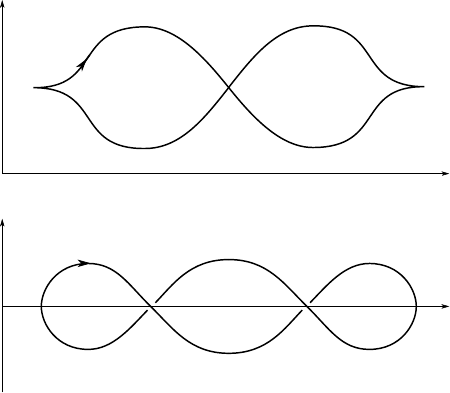}
 \caption{The Legendrian unknot $\Lambda_{-1,-2}$ of ${\tt rot}=-1$ and ${\tt tb}=-2$ depicted in a contact Darboux ball $(\R^2_{xy} \times \R_z,dz-ydx)$. On top we see the front projection, and on the bottom the so-called Lagrangian projection.}
 \label{fig:legendrian}
\end{figure}
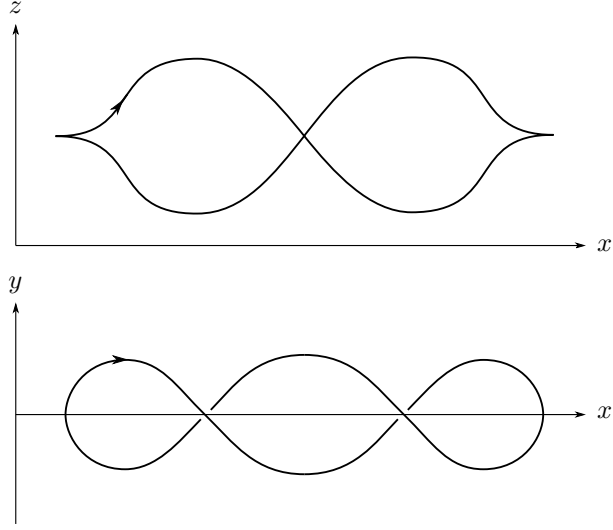

We begin with a local description of the hyperbolic tangency from a symplectic viewpoint. Our description is close to that of Nemirovski--Siegel in \cite[Section 4.3]{NemirovskiSiegel}, where it was shown that the Legendrian unknot $\Lambda_{\pm 1,-2} \subset (S^3,\alpha_{st})$ of ${\tt rot}=\pm1$ and ${\tt tb}=-2$ bounds a disc with a single hyperbolic complex tangency inside the standard ball $B^4 \subset \C^2$. This Legendrian unknot has a representative with the front projection depicted in Figure \ref{fig:legendrian}.

Next we give a description of the local model in terms of the coordinates induced by the standard momentum map
\begin{gather*}
 \C^2 \to (\R_{\ge 0})^2,\\
 (z_1,z_2)\mapsto(\|z_1\|^2,\|z_2\|^2).
\end{gather*}
These coordinates also turn out to be useful for constructing the locally defined plurisubharmonic functions that we will need.

Consider the torus knot
$$T_{-1,2}\coloneqq \left\{(\|z_1\|^2,\|z_2\|^2)=(2/3,1/3), \:\: \OP{arg}(z_2)=-2\OP{arg}(z_1)\right\} \subset S^3$$
which is entirely contained inside the toric fibre over $\{(\|z_1\|^2,\|z_2\|^2)=(2/3,1/3)\},$ and actually smoothly unknotted. Its cone $\R \cdot T_{-1,2} \subset (\C^2,\omega_0)$ is Lagrangian with an isolated singularity at the origin, from which we deduce that $T_{-1,2} \subset (S^3,\alpha_{st})$ is Legendrian. The $(2,1)$-torus knot is smoothly unknotted. However, it is not Legendrian isotopic to the standard Legendrian unknot of $\tt{tb}=-1$, but rather Legendrian isotopic to the once stabilised Legendrian unknot $\Lambda_{\pm 1, -2}$ shown in Figure \ref{fig:legendrian}. This fact can be seen either by an explicit construction or, alternatively, deduced from the classification of Legendrian unknots \cite{TopTrivLeg} by Eliashberg--Fraser. Recall that the Legendrian isotopy class of unknots in the standard contact sphere are completely classified by its rotation number (which is $\pm 1$) together with Thurston--Bennequin invariant (which is $-2$).

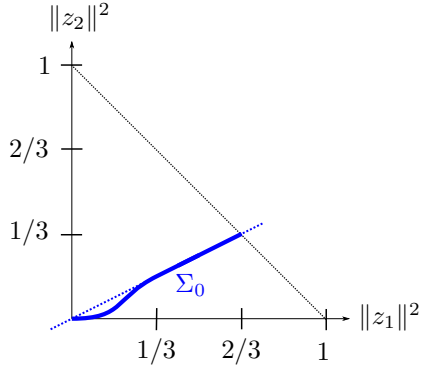
\begin{figure}[htp]
\centering
\vspace{5mm}
\labellist
\pinlabel $1$ at 103 -8
\pinlabel $\|z_1\|^2$ at 128 5
\pinlabel $\|z_2\|^2$ at 11 118
\pinlabel $2/3$ at 72 -8
\pinlabel $1/3$ at -8 37
\pinlabel $2/3$ at -8 68
\pinlabel $1$ at -2 100
\pinlabel $1/3$ at 40 -8
\pinlabel $\color{blue}\Sigma_0$ at 53 17
\endlabellist
\includegraphics[scale=1]{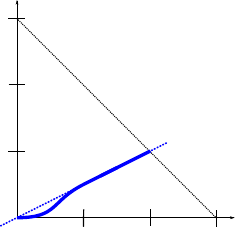}
\vspace{3mm}
\caption{The image of the local model $\Sigma_0$ under the standard momentum map. The projection of $\Sigma_0$ coincides with the line with slope $1/2$ for $\|z_1\|^2 \ge 1/3$, while it is tangent to the horizontal axis precisely at the origin.}
\label{fig:toric}
\end{figure}

It will be useful to introduce some notation for describing surfaces in the Darboux ball in terms of the standard momentum map. First we consider the image of the Lagrangian cone
$$\R \cdot T_{-1,2} \subset (\C^2,\omega_0),$$
with an isolated singularity at the origin, under the standard momentum map $\C^2 \to \left(\R_{\ge 0}\right)^2$; it is the dashed line with slope $1/2$ depicted in Figure \ref{fig:toric}. Note that the cone is not the full pre-image of the line; each fibre consists of the linear sub-torus $\{\OP{arg}(z_2)=-2\OP{arg}(z_1)\}$ that corresponds to the \emph{normal} of the line in the momentum polytope.

We then consider general, not necessarily Lagrangian, surfaces of the form
$$\Sigma(f) \coloneqq \left\{ (z_1,z_2) \in \C^2; \:\|z_2\|^2=f(\|z_1\|^2)\right\} \cap \{\OP{arg}(z_2)=-2\OP{arg}(z_1)\}$$
for some function $f \colon \R_{\ge0} \to \R_{\ge0}$ and note that the aforementioned cone can be described as
$$ \Sigma\left(\frac{1}{2}t\right)=\R \cdot T_{-1,2}$$
defined by the function $f(\|z_1\|^2)=\frac{1}{2}\|z_1\|^2.$

 The main use of this construction will be for creating a particular smoothing of the Lagrangian cone with a singularity at the origin, which is a smooth disc that is totally real except for an isolated hyperbolic complex tangency at the origin, and which satisfies the assumptions of Theorem \ref{thm:rationalconvexity} near the origin.

\begin{lma}
 \label{lma:conesmoothing}
 The Lagrangian cone
 $$\Sigma\left(\frac{1}{2}t\right)=\R \cdot T_{-1,2} \subset (\C^2,\omega_0),$$
which has a singularity precisely at the origin, can be deformed in an arbitrarily small neighborhood of the origin to smooth disc of the form $\Sigma_0=\Sigma(f_0),$
where
 \begin{itemize}
\item $f_{ 0}(t)=t^2$ near $t=0$;
\item $f_{ 0}(t)=\frac{1}{2}t$ for $t \ge 1/3$;
\item $f_{ 0}'(t)> 0$ for all $t > 0$; and
\item $f_{ 0}(t) \le t$ for all $t > 0$.
\end{itemize}
Moreover, this disc is totally real everywhere except near the origin, where it is biholomorphic to a subset of the origin of the graphical surface
$$\left\{v=u^2+\overline{u^2}\right\}=\left\{v=2\mathfrak{Re}(u^2)\right\}=\left\{v=2((\mathfrak{Re}u)^2-(\mathfrak{Im}u)^2)\right\}.$$
In other words, $\Sigma_0$ has a standard hyperbolic complex tangency there.
 \end{lma}
 
\begin{proof}
We just need to verify that $\Sigma_0$ has a standard hyperbolic complex tangency at the origin. To see this, we begin with the straightforward identification of
$$\Sigma (t^2)=\{ (z_1,z_2) \in \C^2; \: \|z_2\|=\|z_1\|^2\} \cap \{\OP{arg}(z_2)=-2\OP{arg}(z_1)\} $$
and the graph
$$z_2=\|z_1\|^4z_1^{-2}= z_1^2\overline{z}_1^2z_1^{-2}= \overline{z_1^2}$$
near the origin. After the coordinate change $(u,v)=(z_1,z_2+z_1^2)$, this graph can be expressed as
$$\left\{v=u^2+\overline{u^2}\right\}=\left\{v=2\mathfrak{Re}(u^2)\right\}=\left\{v=2((\mathfrak{Re}u)^2-(\mathfrak{Im}u)^2)\right\}.$$
This is the so-called standard quadratic hyperbolic tangency, which in particular is ``flat'' in the sense of \cite{Slapar}.
\end{proof}

The next step is to construct a smooth plurisubharmonic function $\rho \colon B^4 \to \R_{\ge0}$ which is strictly plurisubharmonic outside of the origin, and for which $d^c\rho|_{T\Sigma_0}$ is closed, which means that $\Sigma_0$ is Lagrangian inside the subset $B^4 \setminus \{0\}$. In order to do this we first need to describe some useful tools for constructing such a function. We will use polar coordinates $z_i=r_ie^{i\theta_i}$, describe the function through $\rho(t_1,t_2)$ with $t_i=\frac{1}{2}r_i^2$. We compute
$$ -d^c \rho = -d \rho \circ J_0 = r_1^2(\partial_{t_1}\rho(r_1^2/2,r_2^2/2))d\theta_1+r_2^2(\partial_{t_2}\rho(r_1^2/2,r_2^2/2))d\theta_2.$$
Note that, in these coordinates, we can write $\Sigma_0=\{r_2^2=f(r_1^2),\:\:\theta_2=-2\theta_1\}$ for $f$ as described above.

\begin{itemize}
\item \emph{Making $\Sigma_0$ Lagrangian:} The condition for $d^c \rho$ to vanish on $T\Sigma_0$ is that the gradient $\nabla\rho(t_1,t_2)=(\partial_{t_1}\rho,\partial_{t_2}\rho)$ is orthogonal to
 $$2(t_1,-2t_2)=(r_1^2,-2r_2^2)=(r_1^2,-2f(r_1^2))$$ along the subset $\left\{r_2^2=f(r_1^2),\: \theta_2=-2\theta_1\right\}$ or, equivalently, that the gradient is colinear with
 $$2(2t_2,t_1)=(2f(r_1^2),r_1^2),$$
 along the same subset. (The curve $(r_1^2,f(r_1^2))$ shown in Figure \ref{fig:toric} is the projection of $\Sigma_0 \subset B^4$ under the momentum map $(z_1,z_2) \mapsto (r_1^2,r_2^2)$).

\item \emph{Making $\rho$ plurisubharmonic:} Since $\partial_{t_i}=r_i^{-1}\partial_{r_i}$, a direct calculation implies that the function $\rho$ is strictly plurisubharmonic at all points where the inequalities
 $$\frac{1}{r_i}(\partial_{r_i}(r_i^2\partial_{t_i}\rho))=2\partial_{t_i}\rho+r_i^2\partial_{t_i}^2\rho>0$$
 are satisfied simultaneously for $i=1,2.$ In particular, this is automatically the case when $\partial_t \rho>0$ as well as $\partial_{t_i}^2 \rho \ge 0$ are satisfied for $i=1,2$.
 
Furthermore, at a point $p$ where $\partial_{t_i}\rho > 0$ holds simultaneously for both $i=1,2$, strict plurisubharmonicity can be achieved by post-composing $\rho$ with a convex function $g \colon \R \to \R$ for which $g \ge 0$, $g',g''>0$, and where $\frac{d}{dt}\log(g')=g''/g' \gg 0$ is sufficiently large at $t=\rho(p)$. Indeed, the second derivative of $g \circ \rho(t_1,t_2)$ satisfies
 $$\partial_{t_i}^2(g\circ \rho)=g'(\rho)\cdot\partial_{t_i}^2\rho +g''(\rho)\cdot (\partial_{t_i}\rho)^2.$$
\item \emph{Condition for obtaining the standard symplectic form:} When $\rho(t_1,t_2)=C(t_1+t_2)$, it is clear that $g \circ \rho(r_1^2/2,r_2^2/2)$ is strictly plurisubharmonic outside of the origin whenever $g'' \ge 0$ and $g,g' > 0$ hold in the same subset. Recall that the potential
 $$\rho=\rho_0=\frac{1}{2}\left(r_1^2+r_2^2\right)$$
 gives the standard symplectic form $i\del\delbar\rho_0= -dd^c \frac{\rho_0}{2} =\omega_0$.
 \end{itemize}

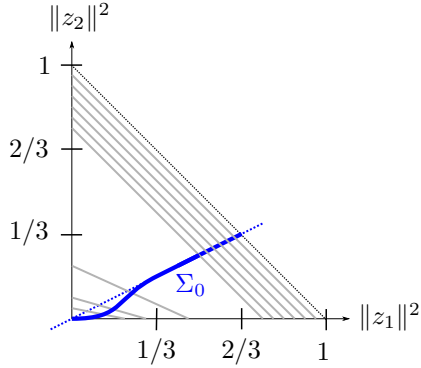
\begin{figure}[htp]
\centering
\vspace{5mm}
\labellist
\pinlabel $1$ at 103 -8
\pinlabel $\|z_1\|^2$ at 128 5
\pinlabel $\|z_2\|^2$ at 11 118
\pinlabel $2/3$ at 72 -8
\pinlabel $1/3$ at -8 37
\pinlabel $2/3$ at -8 68
\pinlabel $1$ at -2 100
\pinlabel $1/3$ at 40 -8
\pinlabel $\color{blue}\Sigma_0$ at 53 17
\endlabellist
\includegraphics[scale=1]{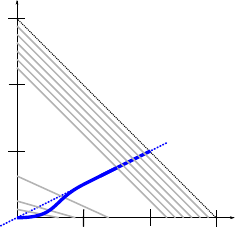}
\vspace{3mm}
\caption{A foliation by convex spheres whose normals are colinear with $(2f(r_1^2),r_1^2)$ along $\Sigma_0=\{r_2^2=f(r_1^2)\}$, and which coincide with the standard concentric round spheres near the boundary of the unit ball. Note that the image of $\Sigma_0$ is the graph of the parabola $\{r_2^2=(r_1^2)^2\}$ near $r_1^2=0$.}
\label{fig:torichypersurfaces}
\end{figure}

The following technical result will be important for establishing the conditions in \cite[Lemme 1]{Gayet}, which is a result that we invoke in the proof of Theorem \ref{thm:rationalconvexity} to show rational convexity.
\begin{prp}
 \label{prp:hyperbolicpsh}
 Let $\Sigma_0 \subset B^4$ be the smooth disc constructed in Lemma \ref{lma:conesmoothing}. There exists a smooth plurisubharmonic function $\rho \colon B^4 \to \R_{\ge 0}$ that satisfies
 \begin{enumerate}
 \item $\rho(\mathbf{z})=\frac{1}{2}\|\mathbf{z}\|^2+C$ near $\partial B^4$ for some constant $C \ge 0$; 
 \item $\rho$ is strictly plurisubharmonic in $B^4 \setminus \{0\}$;
 \item $\rho$ has no critical points in $B^4 \setminus \{0\}$, and its level sets are transverse to $\Sigma_0$ in the same subset; and
 \item the one-form $-\frac{1}{2}d^c \rho$, i.e.~a primitive of $i\del\delbar\rho$ for the exterior derivative, vanishes on $T\Sigma_0$. In particular, $\Sigma_0 \setminus \{0\} \subset (B^4 \setminus \{0\},i\partial\overline{\partial}\rho)$ is an exact Lagrangian submanifold that, for each $c > 0$, intersects the contact hypersurfaces
 $$\left(\rho^{-1}(c),-\frac{1}{2}d^c \rho|_{T\rho^{-1}(c)}\right)$$
 in Legendrian knots.
 \end{enumerate}
 \end{prp}
 \begin{proof}
We start by constructing a smooth function $\rho(t_1,t_2)$ defined on the first quadrant $(\R_{\ge 0})^2$ of the $(t_1,t_2)$-plane whose gradient $\nabla \rho(t_1,t_2)$ is orthogonal to $(t_1,-2f(t_1))$ along the curve $\{t_2=f(t_1)\}$, and such that $\partial_{t_i}\rho>0$ is satisfied away from the origin $(t_1,t_2)=0$ for the partial derivatives of both coordinates $t_i$ with $i=1,2$. 

First consider the foliation of the first quadrant by the lines
$$\left\{s \mapsto ((1+s)t,(1-2s)f(t)), \: s \in [-1,1/2]\right\}_{t \ge 0},$$
where $t$ parametrizes the leaf space (i.e.~the space of lines), and $s$ parametrizes each leaf (i.e.~line). Note that the line with parameter value $t$ passes through the point $(t,f(t))$ when $s=0$ with tangent vector $(t,-2f(t))$.

The line corresponding to the parameter-value $t$ intersects the first and second coordinate axis at the points $(3t/2,0)$ (when $s=1/2$) and $(0,3f(t))$ (when $s=-1$), respectively. Recall that $f'(t) \ge 0$ holds with strict inequality for $t > 0$. Since the $t$-derivative of this family of lines $(1+s,(1-2s)f'(t))$ thus is everywhere non-zero when $s \in [-1,1/2]$ and $t > 0$, we conclude that this is a smooth foliation away from the origin.

The foliation property allows us to find a continuous function $\rho(t_1,t_2)$ defined in the first quadrant $(\R_{\ge0})^2$, which is smooth away from the origin, and whose level-sets constitute the above foliation by lines. We define $\rho$ uniquely by the requirement that it takes the value $3t/2$ along the leaf with parameter value $t$.

The sought plurisubharmonic function that is strictly plurisubharmonic away from the origin will be obtained by setting $t_i=r_i^2/2$ and considering the composition $g \circ \rho$ for a suitable function $g \colon \R \to \R_{\ge 0}$ with $g(0)=0$ and $g'(t) >0$ for $t > 0$.

\emph{Property (1):} Since $f(t)=t/2$ when $t \ge 1/3$, we get the equality $\rho(t_1,t_2)=t_1+t_2$ in the subset $\{t_1+t_2 \ge 1/3\}$ with our convention. We thus need $g(t)=t+C$ to hold in the subset $t \ge 1/3$. 

\emph{Property (2):} First we need to investigate the behavior of the partial derivatives of $\rho(t_1,t_2)$ near the origin $(t_1,t_2)=0$, which is precisely where the smoothness fails. Since $f(t)=t^2$ holds near $t=0$, each level set $\rho^{-1}(t)$ can be parametrized by
$$ s\mapsto (s,t^2(1-s\cdot2/(3t)))=(s,t^2-2st/3)=(s,t(t-2s/3))$$
near the origin. It follows that
$$ t = \rho(t_1,t_2)=t_1/3+\sqrt{t_1^2/9+t_2},\:\:t_1=s, \:\:t_2=t(t-2s/3),$$
which is a function that is smooth away from the origin, has partial differentials that can be bounded from above by functions of the form
\begin{equation}
\label{eq:bound}
\frac{\partial^k\rho }{\partial t_{i_1} \partial t_{i_2} \cdots \partial t_{i_k}} \le C_k\rho^{-m_k}
\end{equation}
for $m_k,C_k > 0$ depending on the order of the partial differential.

For any smooth function $g \colon \R_{\ge 0} \to \R_{\ge 0}$ which satisfies $g(0)=0$, $g^{(k)}(s) \le C_ke^{\frac{-1}{s}}s^{-m_k}$ near $s=0$ for some $C_k > 0$, $m_k > 0$, and all $k=1,2,3,\ldots$, one readily computes that $g \circ \rho$ is smooth in all of $(\R_{\ge 0})^2.$ Here we need Inequality \eqref{eq:bound} for bounding the growth of the the partial differentials of $\rho$. A function $g(t)$ satisfying this property can be obtained by specifying
$$g(0)=0 \:\:\text{and}\:\: g'(t)=e^{\int_{1/2}^t G(s) ds},\:\: t \in [0,1/2],$$
where $G(s) \ge 0$ satisfies $G(s)=1/s^h$ for all $s>0$ close to $s=0$, and $h > 0$. Furthermore, we will take the function to satisfy $G(t)=0$ near $t=1/2$. It follows that $g'(s)>0$ for $s >0$, $g^{(k)}(0)=0$ for all $k =0,1,2,3,\ldots$, and that $g(t)$ is everywhere $C^\infty$. Furthermore, we have $g(t)=t+C$ near $t=1/2$. The latter means that $g\circ \rho(t_1,t_2)=t_1+t_2+C > 0 $ is still satisfied near $\{t_1+t_2=1/2\}$.

After choosing $G(t) \ge 0$ to be sufficiently large away from some neighbourhood of $t = 1/2$ which, in particular means that we must choose $h \gg 0$ sufficiently large, we obtain the strict inequalities
\begin{equation}
\label{eq:ineq}
\partial_{t_i}^2\rho(t_1,t_2)+G(\rho)\cdot (\partial_{t_i}\rho)^2 > 0,\:\:i=1,2,
\end{equation}
away from the origin. To that end, recall that the gradient $\nabla \rho(t_1,t_2)$ has positive components away from the origin, and that Inequality \eqref{eq:bound} controls their growth near the origin. Since $G(t)=g''(t)/g'(t)$, when combined with the second bullet point before this proposition, Inequality \eqref{eq:ineq} implies that $g \circ \rho$ is strictly plurisubharmonic away from the origin.

\emph{Property (3):} Note that the gradient $\nabla\rho(t_1,t_2)$ is orthogonal to the tangent vectors of the level-sets of $\rho$, which by construction is the tangent vectors $(t,-2f(t))$ to the family of lines constructed above. Hence, the gradient is co-linear with $(2f(t),t)$. In particular, since $t > 0$ holds away from the origin, both components of the gradient are positive in this subset. It also follows that the level sets are transverse to $\Sigma_0$. The same is then true for the composition $g \circ \rho$.

\emph{Property (4):} The first bullet point in the paragraph before this proposition, together with the calculation of the gradient $\nabla\rho(t_1,t_2)$, ensures that $d^c\rho$ vanishes when pulled back to $\Sigma_0$. The same is then true for the composition $g \circ \rho$.
\end{proof}

 \subsection{Digression: Results about Lagrangian fillings of $\Lambda_{\pm1,-2}$.}
 \label{sec:digression}

For completeness, we here investigate some properties satisfied by the Lagrangian fillings of the stabilised Legendrian unknot $\Lambda_{\pm 1,-2}$. There are implications for the possible smoothings of the singular Lagrangian cone $\R\cdot\Lambda_{\pm1,-2} \subset (\C^2,\omega_0)$ over the Legendrian knot $\Lambda_{\pm1,-2} \subset (S^3,\alpha_{st})$ shown in Figure \ref{fig:legendrian}.

\begin{figure}[htp]
 \vspace{3mm}
 \label{fig:immersedcobordism}
  \includegraphics{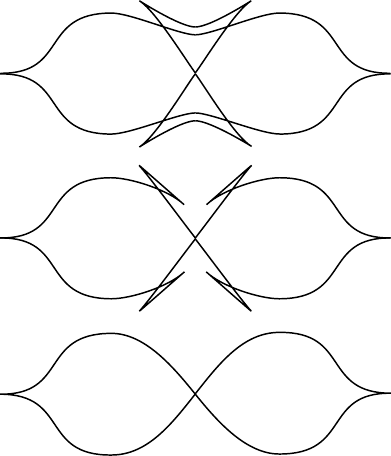}
 \caption{Bottom to middle: A Legendrian isotopy consisting of four Reidemeister one-moves in the front. Middle to top: two Legendrian ambient surgeries performed at the pairs of cusp-edges that face each other produces a link of two standard Legendrian unknots shown on top. This link is a Legendrian realisation of the Whitehead link, and bounds an immersed Lagrangian filling consisting of two discs that intersect transversely in two points with opposite signs.} 
\end{figure}

 \begin{enumerate}
 \item The Legendrian $\Lambda_{\pm1,-2}$ does not bound any orientable Lagrangian submanifolds simply for the purely topological reason that $\tt rot = \pm 1 \neq 0$. Hence, the singularity at the origin of the K\"{a}hler form produced by Proposition \ref{prp:hyperbolicpsh} cannot be avoided. Without a singularity, we would have produced a Lagrangian disc filling of $\Lambda_{\pm1,-2}$ inside the symplectic ball.
 \item The Legendrian $\Lambda_{\pm1,-2}$ also does not bound any non-orientable exact Lagrangian fillings (recall that exact means that $\omega_0$ has a primitive $\eta$ which pulls back to an exact one-form on the Lagrangian). Unlike (1) above, the proof of this fact uses harder techniques than just classical homotopy theory. The Legendrian contact homology of the knot $\Lambda_{\pm1,-2}$, i.e.~the Legendrian isotopy invariant in the form of a differential graded algebra (DGA) defined by Chekanov in \cite{DiffAlg}, is acyclic with $\Z_2$-coefficients. An exact Lagrangian filling would imply the existence of an augmentation of this DGA, i.e.~a unital DGA-morphism to the ground field $\Z_2$; see e.g.~\cite{Chen:NonOrientable}.
 \item The Legendrian knot $\Lambda_{\pm1,-2}$ does, however, bound non-orientable Lagrangian fillings by the h-principle for Lagrangian immersions. Namely, the h-principle provides Lagrangian immersions with generic double points. The latter double points can be removed by a Lagrangian surgery \cite{Polterovich:Surgery}. Recall that performing surgery on several double points typically produces a non-orientable Lagrangian.

 A more explicit construction can be made as follows. Perform ambient Legendrian surgeries on $\Lambda_{\pm1,-2}$ as defined in e.g.~\cite{Dimitroglou:Ambient} to yield an exact Lagrangian cobordism to a link of standard unknots; see Figure \ref{fig:immersedcobordism}. This cobordism has $\Lambda_{\pm1,-2}$ at its concave end and the link of unknots at its convex end. The direction of this cobordism can be reversed, with the cost of adding double points. The reversed immersed cobordism now has two Legendrian unknot components at the concave end; they can be filled by two Lagrangian discs, which intersect since the unknots are linked. Finally, we smooth all double points on the immersed Lagrangian cobordism, which yields a non-exact non-oriented Lagrangian filling.
 \item The Legendrian $\Lambda_{\pm1,-2}$ does not bound a Lagrangian M\"{o}bius band. The reason is that the primitive $-\frac{1}{4}d^c\|\mathbf{z}\|^2$ of the symplectic form, which is exact near the Legendrian boundary, would automatically be exact on the entire Lagrangian M\"{o}bius band; namely, the restriction of $H^1$-cohomology classes to the boundary is an isomorphism when $\R$-coefficients are used. Or put differently: any Lagrangian M\"{o}bius band is automatically exact, which is not allowed by (2) above.

 However, there exists Lagrangian M\"{o}bius bands which are non-exact and non-compactly supported Lagrangian deformations of the cone over $\Lambda_{\pm1,-2}$. In other words, this is a non-exact Lagrangian filling that is not cylindrical, but merely asymptotic to the Legendrian knot at infinity. Such a Lagrangian M\"{o}bius band is given by
 $$ \left\{ \|z_2\|^{2}=\|z_1\|^{-4}+\epsilon\right\} \cap \left\{\OP{arg}(z_2) =-2\OP{arg}(z_1)\right\}$$
 for $\epsilon>0$, whose core is the circle $\{0\}\times\sqrt{\epsilon}\cdot S^1$ contained in the second coordinate plane. (In particular, there is a holomorphic Maslov-0 disc with boundary on this filling.) Note that this Lagrangian is asymptotic to the cone over $\Lambda_{\pm1,-2}$, but that it is neither conical nor exact outside of a compact subset. 
 \end{enumerate}

\begin{rmk} The non-exact Lagrangian asymptotic fillings of the Legendrian knot $\Lambda_{\pm1,-2}$ in (4) above can be used to deform the genus-$g$ surfaces produced by Theorem \ref{thm:exactsurfaces}, which are Lagrangian outside a $2(g-1)$ number of standard models of hyperbolic points, to yield an embedded Lagrangian surface with the hyperbolic points replaced by a $2(g-1)$ number of Lagrangian M\"{o}bius bands. However, since the M\"{o}bius band constructed in (4) does not coincide with the cone over the Legendrian near its boundary -- rather, it is a non-exact deformation of this cone -- we must perform a non-local perturbation of the initial Lagrangian in order for the Lagrangian M\"{o}bius bands to fit in the model. The necessary deformation of the Lagrangian surface can readily be carried out by the fibre-wise addition of a suitable non-exact closed one-form inside its standard Weinstein neighborhood.
\end{rmk}

\section{Non-fillable rationally convex surfaces with only hyperbolic tangencies}
\label{sec:nonfillable}

It is a standard fact that a closed surface $\Sigma \subset \C^2$ of genus $g$ that is totally real except for $h$ hyperbolic tangencies satisfies
$$ \chi(\Sigma)=2-2g=-h.$$
This can e.g.~be seen by computing the self-linking numbers of the knots defined by intersecting the complex tangencies by small round spheres for the framing defined by the Reeb vector field on the round contact sphere $(S^3,\alpha_{st})$, where $\alpha_{st}=\left.\left(-\frac{1}{4}d^c\|\mathbf{z}\|^2\right)\right|_{TS^3}$. In this section we construct non-fillable rationally convex surface in $\C^2$ of genus $g\ge 2$ with a number $2(g-1)>0$ of hyperbolic complex tangencies.

\begin{figure}[htp]
 \vspace{3mm}
 	\labellist
 	\pinlabel $\Lambda_{-2,-2}$ at 30 20
  \pinlabel $\Lambda_{+2,-2}$ at 100 20
  \pinlabel $\Lambda_{0,-3}$ at 65 98
  \pinlabel $\tau$ at 158 145
  \pinlabel $1$ at 170 83
  \pinlabel $\epsilon$ at 170 41
  \pinlabel $0$ at 170 26
	\endlabellist
  \includegraphics[scale=1.2]{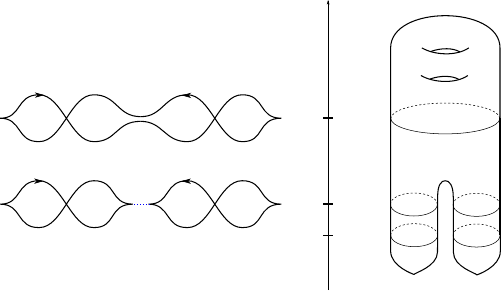}
  \caption{The exact Lagrangian surface with conical singularities produced in the proof of Theorem \ref{thm:exactsurfaces}. In $\{\tau \le 0\}$ we see two Lagrangian cones over the Legendrian unknots $\Lambda_{\pm1,-2}$ contained in disjoint Darboux balls. In $\{\tau < 0\}$ these Darboux balls are connected by a Weinstein handle attachment, and the Lagrangian is a trivial cylinder contained inside of it. A Legendrian ambient surgery produced at the corresponding knots $\Lambda_{\pm1,-2}$ in the standard contact sphere $\{\tau=\epsilon\}$ produces the Legendrian $\Lambda_{0,-3}$ in the contact sphere $\{ \tau=1\}$; there is a Lagrangian pair of pants inside the trivial symplectic cobordism. The top part of the cobordism inside $\{\tau \ge 1\}$ is the exact Lagrangian genus-2 cap constructed by Lin in \cite{Lin}.}
 \label{fig:cobordism}
\end{figure}

The construction of these surfaces can roughly be outlined as follows. We start with an exact Lagrangian genus-$g$ surface inside $(\C^2,\omega_0)$ that is embedded away from a number $2(g-1)$ of conical singularities of the form described in Subsection \ref{sec:model}. Recall that the singularity considered in the aforementioned section is a cone over the Legendrian unknot $\Lambda_{\pm 1,-2} \subset \partial D^4$ of ${\tt tb}=-2$ shown in Figure \ref{fig:legendrian}. It should be noted that the existence of \emph{exact} Lagrangians with all singularities of this particular type is non-trivial; they were first constructed by Lin in \cite{Lin}, whose construction is related to earlier work by Sauvaget \cite{Sauvaget}. Then we replace each cone-point with a smooth disc which has a unique hyperbolic complex tangency, a construction which also was described in Subsection \ref{sec:model}. The concatenation of the Lagrangian and the disc with hyperbolic complex tangencies is depicted schematically in Figure \ref{fig:cobordism}. The construction yields the following:
\begin{thm}
 \label{thm:exactsurfaces}
There exists smooth surfaces $\Sigma_g \subset \C^2$ of any genus $g \ge 2$ that satisfy:
\begin{itemize}
\item There exists a union a $2(g-1)$ number of disjoint round Euclidean K\"{a}hler balls (i.e.~equipped with the standard complex structure $J_0$ and the standard K\"{a}hler form $(i/2)\del\delbar\|\mathbf{z}\|^2$) 
$$\Omega=\underbrace{\left(D^4, J_0 ,(i/2)\del\delbar\|\mathbf{z}\|^2\right) \sqcup \ldots \sqcup \left(D^4, J_0 ,(i/2)\del\delbar\|\mathbf{z}\|^2\right)}_{2(g-1)} \hookrightarrow (\C^2, J_0 ,\omega_0),$$
where $\Sigma_g$ coincides with the standard disc model from Subsection \ref{sec:model} inside each ball. (In particular, in each ball $\Sigma_g$ has a single flat hyperbolic complex tangency, and it is Lagrangian for the K\"{a}hler form $i\del\delbar\rho$ that is degenerate precisely at the tangency); and
\item using $\omega$ to denote the smooth $(1,1)$-form on $\C^2$ that coincides with $i\del\delbar \rho$ in each Darboux ball in $\Omega$ and with $\omega_0$ in $\C^2 \setminus \Omega$, then $\eta|_{T\Sigma_g}$ is exact for any choice of primitive $\eta$ of $\omega=d\eta$. In particular $\Sigma_g$ is an exact Lagrangian submanifold away from the hyperbolic points.
\end{itemize}
 \end{thm}
 \begin{proof}
 Start with a number $2(g-1)$ of disjoint copies of the round closed Darboux balls $(D^4, J_0 ,\omega_0)$, e.g.~obtained by applying affine transformations to one standard Euclidean ball, and denote by $\Omega$ their union. The components of $\Omega$ will be denoted by $D_1,\ldots,D_j,\ldots,D_{2(g-1)} \hookrightarrow \C^2$.
 
 In each boundary $\partial D^4_i$, which is a round contact sphere, we place one copy of the unknot $\Lambda_{1,-2}$ of ${\tt rot}=1$ and ${\tt tb}=-2$. See Figure \ref{fig:legendrian} for the front projection of this knot when placed inside a small contact Darboux ball by a contact isotopy. However, we will use the representative of its Legendrian isotopy class as described in Subsection \ref{sec:model}, i.e.~given as a sub-torus in the standard momentum polytope. To see that these knots are the same up to Legendrian isotopy, one can e.g.~use the fact that unknots inside the standard contact sphere are classified by their classical invariants $\tt{rot}$ and $\tt{tb}$; see \cite{TopTrivLeg}. Finally, denote by $\Lambda \subset \partial \Omega$ the union of these Legendrian knots.

Perform a number $2(g-1)-1$ of standard Weinstein 1-handle attachments in order to connect the symplectic Darboux balls $D_j$ to form a connected ball $D$. These handle-attachments form an abstract compact Weinstein cobordism $(W,d\eta)$ from $\partial \Omega$ to $\partial D$. For very thin Weinstein 1-handles the resulting Weinstein cobordism $(W,d\eta)$ admits a symplectic embedding into $(\C^2,\omega_0)$, with concave end coinciding with the boundary $\partial \Omega$ of our original domain $\Omega \subset \C^2$. We refer the reader to \cite{Cieliebak:SteinWeinstein}, but also provide the following rough outline of a construction:

The Weinstein 1-handles are locally determined by the corresponding isotropic core disks, which constitute the skeleton of the cobordism. The core disks are the 1-dimensional isotropic discs $(\Gamma_0,\partial \Gamma_0) \subset (W,\partial \Omega)$ defined as the stable manifolds of the critical points of the Liouville vector field on $(W,d\eta)$. The embedding of $(W,d\eta)$ can then be constructed as follows. First, we embed a number $2(g-1)-1$ of isotropic arcs $\Gamma \subset \C^2 \setminus \OP{int}\Omega$ with boundary transverse to $\partial \Omega$. Second, we identify $\partial \Omega \cup \Gamma_0 
\subset W$ with $\partial \Omega \cup \Gamma \subset \C^2$ by choosing a smooth identification of $\Gamma_0$ with $\Gamma$. Finally, this identification can be extended to a smooth symplectomorphism of a neighbourhood of $\partial \Omega \cup \Gamma_0$ into $\C^2$.

By the previous paragraph the cobordism can be identified with a symplectic embedding $(W,d\eta) \hookrightarrow (\C^2,\omega_0)$, whose convex boundary
$$Y \coloneqq \partial D \subset \partial W=\partial \Omega \sqcup \partial D$$
thus is a contact-type hypersurface in $(\C^2,\omega_0)$; this contact manifold is contactomorphic to the standard contact 3-sphere. (Note that the contact form $\alpha\coloneqq \eta|_{TY}$ that is induced by the primitive $\eta$ of $\omega_0$ typically is not equal to the round contact form, but this will not be important here.) The concave boundary of the cobordism $W \subset \C^2$ is the disjoint union $\partial \Omega \subset \partial W$ of a number $2(g-1)$ of standard round contact spheres.

The endpoints of the unstable manifolds $\Gamma\subset W$, i.e.~the attaching $0$-spheres of the Weinstein 1-handles, can be assumed to be disjoint from the Legendrian knots $\Lambda \subset \partial \Omega \subset \partial W$ contained in the concave boundary of the cobordism $W$. Furthermore, since the complement of a point in the standard contact sphere is contactomorphic to a Darboux ball (see e.g.~\cite{Geiges:Intro}) we can assume that $\Lambda \subset B \subset \partial \Omega$ is contained in a union of contact Darboux balls $B \subset \partial\Omega$ that are disjoint from the attaching spheres. 

Using the Liouville flow defined by $\eta$ in the Weinstein cobordism $(W,d\eta)$ we obtain an exact embedding of the cylinder
$$([-a,0]_\tau \times B,\{-a\} \times B, \{0\} \times B, e^\tau\alpha_B) \hookrightarrow (W,\partial \Omega,Y,\eta).$$
Inside this trivial symplectic cobordism, we can embed the trivial Lagrangian cobordism $[-a,0]\times \Lambda$ with convex cylindrical end $\Lambda_+ \subset (Y,\alpha)$ equal to a union of a number $2(g-1)$ of $\Lambda_{\pm1,-2}$-knots contained in disjoint contact Darboux balls in the contact sphere $Y$. Note that this trivial Lagrangian cobordism is tangent to the Liouville vector field of $\eta$ by construction.

 We then use Lemma \ref{lem:symplectization} to extend the embedding of $Y$ to a symplectic embedding of
 $$([0,+\infty)_\tau \times Y,d(e^\tau\alpha)) \subset (\C^2,\omega_0),$$
 where $\{\tau =0\}$ is the original embedding of $Y$, and where $\{\tau >0\}$ is contained inside the unbounded component of the complement $\C^2 \setminus Y$. Adjoining the above embedding to $W$ produces a symplectic embedding
 $$(\hat{W},d\eta) \hookrightarrow (\C^2,\omega_0)$$ of the completion $\hat{W}$ of $W$.
 
 We are now ready to construct the surface $\Sigma_g \setminus (\C^2,\omega)$ which will be exact Lagrangian outside of the hyperbolic points.
 
\emph{A disc inside $\Omega$:} Here we use the model from Subsection \ref{sec:model} in each component of $\Omega \subset \C^2$. The components are the discs with one complex hyperbolic tangency, which is Lagrangian for $\omega$ away from these tangencies.
 
\emph{A trivial cylinder in $W$:} Here we use the trivial Lagrangian cobordism $[-a,0] \times \Lambda$ that is tangent to the Liouville vector field.
 
 \emph{A Lagrangian handle-attachment inside $\hat{W} \setminus W$:} First we perform a Legendrian isotopy to make the union $\Lambda_+$ of Legendrian unknots equal to the standard representatives shown in Figure \ref{fig:legendrian} horizontally aligned next to each other. There is an associated exact Lagrangian trace cobordism.
 
Then we perform a number $2(g-2)$ of ambient Legendrian surgeries as defined in e.g.~\cite{Dimitroglou:Ambient} on $\Lambda_+$, thus yielding a single component Legendrian unknot $\Lambda_{0,-2(g-1)-1} \subset Y$; see Figure \ref{fig:cobordism} for the case $g=2$. There is an associated exact Lagrangian cobordism that consists of the corresponding standard Lagrangian one-handles. This Lagrangian cobordism lives in
 $$[0,1] \times Y \subset (0,+\infty) \times Y,e^\tau\alpha) \cong (\hat{W} \setminus W,\eta)$$
 and has a concave end consisting of the unlinked unknots $\{0\} \times \Lambda_+ \subset (Y,\alpha)$ and a convex end given by the connected unknot $\Lambda_{0,-2(g-1)-1} \subset \{1\} \times Y$. Note that the Lagrangian cobordism is diffeomorphic to a sphere with $2(g-1)+1$ discs removed. We refer to \cite{Dimitroglou:Ambient} for an explicit description of its construction. 

\emph{Capping off $\Lambda_{0,-2(g-1)-1}$:} The final piece is the most involved part of the construction, and consists of adjoining an exact Lagrangian cobordism in $[1,+\infty) \times Y$ of genus $g$ with concave end $\Lambda_{0,-2(g-1)-1}$. Here we rely on Lin's construction in \cite{Lin}. Since we later want to understand the smooth isotopy class of the filling we proceed to give some details.
 
\begin{figure}[htp]
 \vspace{3mm}
  \includegraphics[scale=1.2]{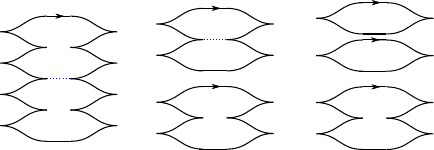}
  \caption{The Legendrian unknot $\Lambda_{0,-2(g-1)-1}$ for $g=4$ (shown on the left) is Lagrangian cobordant to a union of two unlinked Legendrian unknots $\Lambda_{0,-3}$ (shown in the middle) by an exact Lagrangian handle-attachment cobordism that corresponds to the Legendrian ambient surgery along the dashed line on the left. The link shown in the middle is Lagrangian cobordant to a Legendrian link consisting of three unlinked unknots; one copy of $\Lambda_{0,-3}$ together with two standard Legendrian unknots $\Lambda_{0,-1}$. Again, this cobordism is a standard Lagrangian handle-attachment cobordism induced by a Legendrian ambient surgery along the dashed line in the middle.}
\label{fig:stabilisations}
\end{figure}

In the simplest case $g=2$, i.e.~when there are exactly two hyperbolic tangencies, the construction of the exact Lagrangian is depicted in Figure \ref{fig:cobordism}. When $g \ge 2$, we perform a sequence of ambient surgeries along unknotted arcs to produce an exact Lagrangian cobordism of genus $g-2$ from $\Lambda_{0,-2(g-1)-1}$ to $\Lambda_{0,-3}$. This is obtained by first performing a number $2g$ of suitable ambient Legendrian surgeries to produce an exact Lagrangian cobordism from $\Lambda_{0,-2(g-1)-1}$ to a union of unlinked Legendrian unknots: one copy of $\Lambda_{0,-3}$ and $2g$ copies of the standard Legendrian unknot $\Lambda_{0,-1}$; see Figure \ref{fig:stabilisations} for the case $g=4$. Then, we perform ambient Legendrian surgeries from the latter link of Legendrian unknots to $\Lambda_{0,-3}$ by taking connected sum with the standard Legendrian unknots. (Recall that a cusp-connected sum with a standard Legendrian unknot does not change the Legendrian isotopy class of the second knot.) The resulting cobordism has genus $g-2$. Finally, we adjoin Lin's exact Lagrangian cobordism of genus 2 to the Legendrian knot $\Lambda_{0,-3}$.

The exactness of $\Sigma_g$ is immediate. Take any global primitive one-form $\eta \in \Omega^1(\C^2)$ of the two-form $\omega$. Since $\eta$ is closed on the pull-back $\eta|_{T\Sigma_g}$, it is moreover exact on the discs $\Sigma_g \cap \Omega$. Since $
\omega=\omega_0$ on $\C^2 \setminus \Omega$, it follows that $\eta=-\frac{1}{4}d^c\|\mathbf{z}\|^2+\beta$ for some closed one-form $\beta$ on the same domain. Since $\C^2 \setminus \Omega$ is simply connected, we get that $\beta=df$ is exact. Since $-\frac{1}{4}d^c\|\mathbf{z}\|^2|_{T\Sigma_g \setminus \Omega}$ is exact by construction, it then follows that $\eta|_{T\Sigma_g}$ is exact as well.
 \end{proof}
 A direct application of Theorem \ref{thm:rationalconvexity}, where we have used the properties of the model disc neighborhood $\Sigma_0$ of each hyperbolic complex tangency proven in Proposition \ref{prp:hyperbolicpsh} , gives that
\begin{cor}
The surfaces with only hyperbolic complex tangencies as constructed by Theorem \ref{thm:exactsurfaces} are all rationally convex.
\end{cor}

Recall the famous result by Gromov \cite{Gromov:Pseudo}, which implies that no real $n$-dimensional closed sub manifold is exact for a global K\"{a}hler form $\omega=d\eta$ on $\C^n$, i.e.~$\eta$ cannot be exact when pulled back to any $n$-dimensional submanifold. In the case when $\omega$ is not required to be a K\"{a}hler form everywhere, but merely non-negative, the situation becomes very different. The surfaces produced by Theorem \ref{thm:exactsurfaces} are all exact in the following weaker sense:
\begin{dfn}
 An embedded submanifold $\Sigma \subset \C^n$ is said to be \emph{exact relative to $S$ }, or alternatively just \emph{exact}, for a smooth subset $S \subset \C^n$, if the following is satisfied:
 \begin{itemize}
 \item any Riemann surface in $S$ with boundary contained in $S \cap \Sigma$ is constant;
 \item there exists a smooth non-negative $(1,1)$-form $\omega=-dd^c\phi \in \Omega^2(\C^n)$ which is non-degenerate (and thus a K\"{a}hler form) away from $S \subset \C^n$; and
 \item there exists a global primitive $\eta \in \Omega^1(\C^n)$ of $\omega=d\eta$ for which $\eta|_{T\Sigma} \in \Omega^1(\Sigma)$ is exact.
 \end{itemize}
\end{dfn}
 
\begin{rmk}
 \label{rmk:totreal}
 Note that the submanifold $\Sigma \setminus S$ is totally real, since it is Lagrangian for the K\"{a}ler structure $(\C^n \setminus S,\omega)$. However, note that nothing prevents $\Sigma$ from having complex tangencies inside $S$.
 \end{rmk}
 
We deduce the following result. 
\begin{lma}
Every compact Riemann surfaces with boundary attached to submanifold $\Sigma \subset \C^n$ which is exact relative to some subset $S \subset \C^n$ must be constant.
\end{lma}
\begin{proof}
 Consider a continuous map from a Riemann surface $u \colon (C,\partial C) \to (\C^n,\Sigma)$ which is holomorphic away from the boundary. It follows by a classical result of Chirka \cite{Chirka} that $u$ constitutes a rectifiable current, which enables us to use Stokes' theorem in the setting of currents. 
 
 Stokes' theorem for currents implies that
 $$\int_u\omega=\int_{\partial u}\eta=\int_{\partial u}\eta|_{T\Sigma}=\int_{\partial u}df=0.$$
 This means that, in the case when $u$ intersects the subset $\C^n \setminus S$ where the two-form $\omega$ is non-degenerate, the map $u$ must be constant since otherwise $\int_u \omega>0.$ On the other hand, if $u$ is contained entirely inside $S$, then it is constant by the assumption of $\Sigma$ being exact relative $S$.
\end{proof}
 In other words, the surfaces produced by Theorem \ref{thm:exactsurfaces} do not admit holomorphic fillings in the following very strong sense.
 \begin{cor}
 \label{cor:rationallyconvex}
 Any continuous map from a compact Riemann surface into $\C^2$ with boundary mapped to one of the rationally convex surfaces $\Sigma_g \subset \C^2$ produced by Theorem \ref{thm:exactsurfaces} must be constant, under the assumption that the map is holomorphic away from the boundary. 
\end{cor}

We end with an auxiliary lemma about contact hypersurfaces in standard symplectic $(\C^n,\omega_0)$. Recall that a primitive $\eta$ of a symplectic form $\omega$ gives rise to the Liouville vector field $\zeta$ defined by $\iota_\zeta\omega=\eta$.
\begin{lma}
 \label{lem:symplectization}
 If $Y \subset (\C^n,\omega_0)$ is a closed hypersurface that satisfies $H_1(Y;\R)=0$ which is of contact type, i.e.~$\omega_0$ admits a primitive one-form $\eta$ near $Y$ whose Liouville flow is transverse to $Y$, then there exists a symplectic embedding of the so-called symplectization
 $$ (\R_\tau \times Y,d(e^\tau \alpha)) \hookrightarrow (\C^n,\omega_0) $$
 whose restriction to $\{\tau=0\}$ is the original embedding of $Y$, and where $\alpha \coloneqq \eta|_{TY}$ is the contact one-form on $Y$.
\end{lma}
\begin{proof}
 Any two primitives of $\omega_0$ differ by a closed one-form near $Y$. Since $H_1(Y)=0$, this one-form is exact, and we can thus extend $\eta$ to a global primitive of $\omega_0$ on $\C^n$ that coincides with the standard primitive $-\frac{1}{4}d^c\|\mathbf{z}\|^2$ of $\omega_0$ outside of a compact subset. One can then readily use the Liouville flow in order to construct the sought symplectization coordinates. The Liouville flow is complete since the primitive of the symplectic form is given by $-\frac{1}{4}d^c\|\mathbf{z}\|^2$ outside of a compact subset.
 \end{proof}

\section{Constructing knotted rationally convex surfaces}
Here we construct examples of smoothly knotted surfaces in $\C^2$ that are rationally convex and which are exact with only hyperbolic complex tangencies. A closed embedded surface in $\C^2$ that bounds an embedded three-dimensional handlebody is called \textbf{unknotted}. Recall that the smooth isotopy classes of surfaces in $\C^2$ form a very rich structure, that can be detected using the fundamental group of the surface complement. The set of such fundamental groups are as rich as the knot groups of one-dimensional knots in $\R^3$, while the complement of an unknotted surface is equal to $\Z$ by a direct computation; see e.g.~\cite[Lemma 3.1]{Dimitroglou:Approximation}.

First we show that
\begin{prp}
The rationally convex surfaces $\Sigma_g \subset \C^2$ with only hyperbolic tangencies constructed in the proof of Theorem \ref{thm:exactsurfaces} all bound embedded three-dimensional handlebodies; i.e.~they are unknotted.
\end{prp}
\begin{proof}
In order to produce the required handlebody that bounds the surface, i.e.~establish unknottedness, it suffices to show the existence of a proper Morse function $f \colon \C^2 \to \R_{\ge0}$ with a unique critical point -- a local minimum -- for which the following is satisfied:
\begin{itemize}
\item $\Sigma_g$ is disjoint from the critical point of $f$ and the restriction $f|_{\Sigma_g}$ to the surface is Morse;
\item the intersection $\Sigma_g \cap f^{-1}(t) \subset f^{-1}(t) \cong S^3$ is a union of unlinked unknots for all regular values $t$ of $f|_{\Sigma_g}$; 
\end{itemize}
and, most importantly, 
\begin{itemize}
\item the unlinked unknots in the regular levels of $f|_{\Sigma_g}^{-1}(t)$ bound a smoothly varying embedded union of discs inside $f^{-1}(t)$ that become pinched or split exactly at the critical points of $f|_{\Sigma_g}$. 
\end{itemize}
Note that any surface can be put in a position satisfying the first two bullet points; it is only using the last bullet point that we can ensure unknottedness. The mechanism that we will use to ensure that the last bullet point holds is that
\begin{itemize}
\item the regular slices $f|_{\Sigma_g}^{-1}(t) \subset f^{-1}(t)$ have a smooth family of knot projections to $\R^2$ under which two different components have disjoint images, and where each image of a component moreover is contained in the unbounded component of the complement of the others.
\end{itemize}

Next we produce such a Morse function for the surface $\Sigma_g$. Recall that this surface can be decomposed into an unknotted disc contained inside the smooth ball $\Omega \cup W$, together with a genus-$g$ surface contained inside $\C^2 \setminus (\Omega \cup W)$; we refer to the proof of Theorem \ref{thm:exactsurfaces} for the construction of the ball $\Omega \cup W \subset \C^2$.

It is easy to construct the sought Morse function inside the smooth ball $\Omega \cup W$. What remains is to check that this Morse function extends to a proper Morse function to the complement
$$\C^2 \setminus (\Omega \cup W) \cong \R \times S^3$$
of the ball, such that the restriction to $\Sigma_g$ satisfies the last bullet point above.

We claim that the sought Morse function on $\C^2 \setminus (\Omega \cup W)$ can simply be taken to be the symplectisation-coordinate $\tau$. It can be seen explicitly by our construction that the part of the Lagrangian cobordism that we construct using explicit Legendrian ambient surgeries in the proof of Theorem \ref{thm:exactsurfaces} has knot projections with the required properties; see e.g.~the front projections Figures \ref{fig:cobordism} and \ref{fig:stabilisations}.

The top piece of $\Sigma_g$, i.e.~the exact Lagrangian genus-two cap produced by Lin in \cite{Lin} also consists of slices whose knot projections satisfy the sought properties, as can be seen by inspecting \cite[Figures 23, 24 and 25]{Lin} in his construction; these figures depict knot projections of the generic slices of his surface, and they clearly satisfy the property of the last bullet point above.
\end{proof}
We then show that the surface produced by Theorem \ref{thm:exactsurfaces} can be deformed in order to yield rational convex surfaces in different knot classes, under the assumption that $g \gg 0$ is sufficiently large.
\begin{thm}
\label{thm:knotted}
For each fixed $g \gg 0$ sufficiently large, there exist several non-isotopic embeddings of the genus-$g$ surface inside $\C^2$ which all are rationally convex, exact, and have only hyperbolic complex tangencies. Moreover, for any $k > 0$, the fundamental groups of the complements of these different embeddings can be assumed to live in $k$ pairwise different isomorphism classes.
\end{thm}
This result should be contrasted to the case of rationally convex surfaces in $\C^2$ which are totally real and orientable; these are all tori that are smoothly unknotted by work of the first author joint with Goodman and Ivrii \cite{Dimitroglou:Isotopy}.
\begin{cor}
\label{cor:knotted}
There are non-orientable surfaces that admit several Lagrangian embeddings in $\C^2$ that are different when considered up to smooth isotopy. (Recall that these embeddings in particular are totally real and rationally convex by \cite{DuvalSibony}.)
\end{cor}
\begin{proof}[Proof of Corollary \ref{cor:knotted}]
Each hyperbolic point is contained inside a small four-ball, the union of which we denote by $\Omega$, such that the surface outside of the union of four-balls $\Sigma_g \setminus \Omega$ is Lagrangian for the standard symplectic structure $\omega_0$. Furthermore, the intersection $\Sigma_g \cap \Omega$ is a union of unknotted discs that intersect the boundary components of $\Omega$ (i.e.~a disjoint union of standard contact spheres) in Legendrian unknots.

We claim that it is sufficient to replace the unknotted discs in $\Omega$ with non-orientable surfaces that are Lagrangian for the standard symplectic form $\omega_0$, coincide with $\Sigma_g$ near $\partial \Omega$, and which are unknotted in the sense that they bound a handlebody with a half-disk removed. More precisely, we want there to exist a handlebody in $\C^2$ whose boundary intersects $\Omega$ precisely in the Lagrangian surface. In particular, it follows that the complement of the surface inside the ball has fundamental group $\Z$ (which is isomorphic to the fundamental group of the complement of the unknotted disc).

The sought non-orientable Lagrangian surface can readily be constructed using Item (3) in Subsection \ref{sec:digression}. The knottedness then follows by a computation of the fundamental group of the complement of the non-orientable Lagrangian constructed using the Seifert--van Kampen theorem.
\end{proof}
\begin{proof}[Proof of Theorem \ref{thm:knotted}]
Recall that the construction of the rationally convex surface $\Sigma_g \subset \C^2$ in the proof of Theorem \ref{thm:exactsurfaces} is carried out by adjoining an exact Lagrangian cap of genus $g\ge 2$ to a Legendrian unknot $\Lambda_{0,2(g-1)-1} \subset Y \cong S^3$ of ${\tt rot}=0$ and ${\tt tb}=-2(g-1)-1$ with $g-1$ positive and $g-1$ negative stabilisations. More precisely, by construction
$$\Lambda_{0,2(g-1)-1} \subset Y = \partial (\Omega \cup W) \subset \C^2$$
is contained inside a contact-type hypersurface $Y \cong S^3$ in $(\C^2,\omega_0)$ that is the boundary of a 4-ball $\Omega \cup W$; see the proof of Theorem \ref{thm:exactsurfaces}. Moreover, the piece $\Sigma_g \cap (\Omega \cup W)$ of the surface contained in the ball is an unknotted disc with $g-1$ hyperbolic complex tangencies, while the genus-$g$ cap $\Sigma_g \cap (\C^2 \setminus \Omega \cup W)$ is an exact Lagrangian that lives in the complement of the ball.

By construction, we may assume that the surface $\Sigma_g$ coincides with the Lagrangian cylinder
$$[1,B] \times \Lambda_{0,2(g-1)-1} \subset ([1,B]_\tau \times Y,d(e^\tau\alpha)) \hookrightarrow (\C^2,\omega_0)$$
in symplectization coordinates defined near the contact-type hypersurface $(Y,\alpha)$. In fact, here we may assume that $B \gg 0$ is arbitrarily large, since we can perform a suitable rescaling of the symplectic form near $Y$ as well as in the unbounded region of $\C^2 \setminus Y$ by an application of the Liouville flow.

The smooth isotopy class of $\Sigma_g$ will be altered by replacing the trivial concordance $[1,B] \times \Lambda_{0,2(g-1)-1} \subset [1,B] \times Y$ with a knotted Lagrangian concordance that coincides with the original trivial cylinder near its boundary.

By the construction of \cite[Theorem 1.2]{Dimitroglou:Approximation} we can find such knotted Lagrangian concordances
$$ ([1,B] \times S^1, \{A\} \times S^1,\{B\}\times S^1) \hookrightarrow ([1,B] \times Y,\{1\} \times \Lambda_{0,2(g-1)-1}, \{B\} \times \Lambda_{0,2(g-1)-1}, d(e^t\alpha)) $$
whose complement has the fundamental groups of any given finite subset of the set of knot groups for 1-knots in $S^3$, under the assumption that we take $g \gg 0$ to be sufficiently large. The reason for this is that the construction needs sufficiently many stabilisations of both signs of the Legendrian knot involved; in our case $\Lambda_{0,2(g-1)-1}$ has $g-1$ positive and $g-1$ negative stabilizations. Since the original surface bounds a handlebody, the fundamental group of the complement of the new deformed surface is the same as the fundamental group of complement of the concordance, which agrees with the given knot group; see \cite[Corollary 3.4]{Dimitroglou:Approximation}. (The latter result was formulated for tori, but the proof is the same also in the case of higher-genus surfaces that bound handlebodies.)
\end{proof}

\bibliographystyle{alpha}
\bibliography{references}

\end{document}